\documentclass{article}
\usepackage{amsmath,amssymb,amsthm}
\usepackage{graphicx,subfigure,float,url,color}
\usepackage{mathrsfs}
\usepackage[colorlinks=true]{hyperref}
\usepackage{pdfsync}
\usepackage{enumitem}

\def\R{\mathrm{I\kern-0.21emR}}
\def\N{\mathrm{I\kern-0.21emN}}

\renewcommand{\geq}{\geqslant}
\renewcommand{\leq}{\leqslant}

\newtheorem{theorem}{Theorem}

\theoremstyle{definition}\newtheorem{example}{Example}
\theoremstyle{definition}\newtheorem{remark}{Remark}

\title{Turnpike in optimal control and beyond: a survey}

\author{
Emmanuel Tr\'elat\footnote{Sorbonne Universit\'e, Universit\'e Paris Cit\'e, CNRS, Inria, Laboratoire Jacques-Louis Lions, LJLL, F-75005 Paris, France (\texttt{emmanuel.trelat@sorbonne-universite.fr}).}
\and
Enrique Zuazua
\footnote{
    [1] Friedrich-Alexander-Universit\"at Erlangen-N\"urnberg, 
    Department of Mathematics, Chair for Dynamics, Control, Machine Learning and Numerics 
    (Alexander von Humboldt Professorship), 
    Cauerstr. 11, 91058 Erlangen, Germany (\texttt{enrique.zuazua@fau.de}).
    \newline\indent
    [2] Chair of Computational Mathematics,
 Deusto University, 48007 Bilbao, Basque Country, Spain.
    \newline\indent
    [3] Universidad Aut\'onoma de Madrid,
    Departamento de Matem\'aticas, 
    Ciudad Universitaria de Cantoblanco, 28049 Madrid, Spain.
}
}

\date{}

\begin{document}

\maketitle

\begin{abstract}
The turnpike principle is a fundamental concept in optimal control theory, stating that for a wide class of long-horizon optimal control problems, the optimal trajectory spends most of its time near a steady-state solution (the ``turnpike") rather than being influenced by the initial or final conditions.

In this article, we provide a survey on the turnpike property in optimal control, adding several recent and novel considerations.
After some historical insights, we present an elementary proof of the exponential turnpike property for linear-quadratic optimal control problems in finite dimension. Next, we show an extension to nonlinear optimal control problems, with a local exponential turnpike property. On simple but meaningful examples, we illustrate the local and global aspects of the turnpike theory, clarifying the global picture and raising new questions. 
We discuss key generalizations, in infinite dimension and other various settings, and review several applications of the turnpike theory across different fields.
\end{abstract}

\tableofcontents

\section{A bit of history: discovery of the turnpike phenomenon}
The turnpike property refers to the often encountered phenomenon that solutions in large time of an optimal control problem tend to spend most of the time near a steady-state (or a more general set) called the turnpike. More generally, it is a behavior widely shared by a number of dynamical optimization problems, where the turnpike, when it is a steady-state, appears to be the solution of a static optimization problem. 

\begin{figure}[h]
\centerline{\includegraphics[width=9cm]{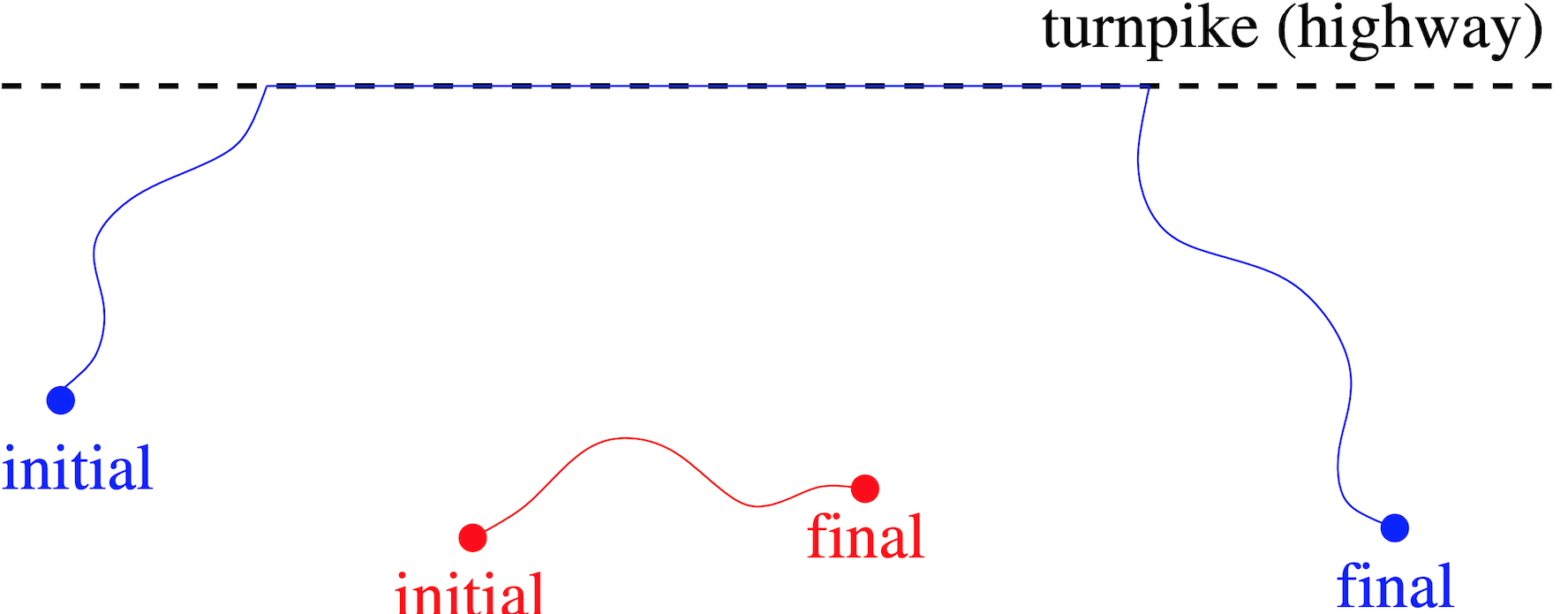}}
\caption{Illustration of the turnpike phenomenon, as described by Samuelson {\it et al.} While the blue trajectory passes through the turnpike, spending most of the time on it, because the initial and final states are far enough, the optimal control strategy in red evolves directly from the initial state towards the final one, because they are close enough.}\label{fig_turnpike}
\end{figure}

\medskip
The turnpike property, at least for the state (not for the costate or adjoint state) has been discovered long time ago. The name ``turnpike" has been coined by the Nobel Prize Paul Samuelson and his co-authors in the article \cite[Chapter 12]{DorfmanSamuelsonSolow} in econometrics, with the objective of deriving efficient programs of capital accumulation, in the context of a Von Neumann model in which labor is treated as an intermediate product. As quoted by \cite{MacKenzie1976}, in \cite[Chapter 12]{DorfmanSamuelsonSolow} one can find the following seminal explanation, illustrated on Figure \ref{fig_turnpike}:
\begin{quote}
\textit{Thus in this unexpected way, we have found a real normative significance for steady growth -- not steady growth in general, but maximal von Neumann growth. It is, in a sense, the single most effective way for the system to grow, so that if we are planning long-run growth, no matter where we start, and where we desire to end up, it will pay in the intermediate stages to get into a growth phase of this kind. It is exactly like a turnpike paralleled by a network of minor roads. There is a fastest route between any two points; and if the origin and destination are close together and far from the turnpike, the best route may not touch the turnpike. But if origin and destination are far enough apart, it will always pay to get on to the turnpike and cover distance at the best rate of travel, even if this means adding a little mileage at either end. The best intermediate capital configuration is one which will grow most rapidly, even if it is not the desired one, it is temporarily optimal.}
\end{quote}
It appears that, prior to this discovery, the phenomenon had already been observed by von Neumann~\cite{vonNeumann} and Ramsey~\cite{Ramsey}. Consequently, in econometrics, the term \emph{von Neumann point} is sometimes used to refer to the \emph{turnpike point} $\bar{x}$.

We also refer to~\cite{Cass1965} for a specific optimal economic growth problem in which this phenomenon was highlighted. 
It was observed that, in many contexts, the solution to an optimal control problem over a large time horizon tends to remain close to a steady state for most of the time. 
In the infinite-horizon case, the solution is expected to converge to this steady state. 
In this setting, the \emph{turnpike property} refers to the fact that optimal trajectories over long time horizons are expected to converge, in an appropriate sense, to the so-called \emph{von Neumann points}.

\medskip
Some turnpike theorems have been derived in the 60's for discrete-time optimal control problems arising in econometrics (see, e.g., \cite{MacKenzie1963}). Continuous versions have been proved in \cite{Haurie} under quite restrictive assumptions on the dynamics motivated by economic growth models. All of them are established for point-to-point optimal control problems and give information on the trajectory only (but not on the adjoint vector).
We also refer the reader to \cite{CarlsonHaurieLeizarowitz_book1991} for an extensive overview of these continuous turnpike results (see also \cite{Zaslavski_2006, Zaslavski_2015}).

\medskip
There has been an impressive revival of the interest to turnpike in the years 2010. The turnpike property has been remarkably used to establish convergence of the Model Predictive Method (MPC) in \cite{Grunediscrete},  and its random-batch version \cite{veldman}. 
Turnpike phenomena have been put in evidence in optimal control problems coming from biology, such as in \cite{CoronGabrielShang, DjemaGiraldiMaslovskayaBernard}, in relation with singular arcs (see also \cite{Rapaport}), in problems related with human locomotion (see \cite{CJM}), with fluid motions (see \cite{Zamorano}) or with ecosystems (see \cite{CaillauDjemaGouzeMaslovskayaPomet}), 
in sports models (see \cite{AftalionTrelat_RSOS, AftalionTrelat_JOMB}), in dynamic optimal design (see \cite{AM, LanceTrelatZuazua_SCL2020}), in mean field games (see \cite{CirantPorretta_COCV2021, Porretta_2018, SunYong_SICON2024}) and even in deep learning (see \cite{EsteveGeshkovskiPighinZuazua_2020, GeshkovskiZuazua_AN2022}),
where optimal trajectories in long time enjoy the above asymptotic property called turnpike phenomenon.
%Note that, in \cite{BoscainPiccoli}, the word `turnpike" refers to the set of points where singular trajectories can stay. In dimension $2$ for the minimal time problem for a control-affine system $\dot x(t)=f_0(x(t))+u(t)f_1(x(t))$, with $u(t)\in[-1,1]$, it is the set of points where $f_1$ is parallel to the Lie bracket $[f_0,f_1]$. This is something different. 

\medskip
Actually, as we will see in the proof of Theorem \ref{thm_turnpike_LQ} in the LQ setting, and as fully exploited in \cite{TrelatZuazua_JDE2015}, the turnpike property is due to the saddle point feature of the extremal equations of optimal control (see \cite{Rockafellar1973, Samuelson1972}), and more precisely to the Hamiltonian nature of the extremal equations inferred from the Pontryagin maximum. These results relate the turnpike property with the asymptotic stability properties of the solutions of the Hamiltonian extremal system, coming from the concavity-convexity of the Hamiltonian function.
These facts were pointed out in the articles \cite{AndersonKokotovic,WildeKokotovic}, which do not seem to have been widely acknowledged. Yet, the authors of these papers proved that the optimal trajectory is approximately made of two solutions for two infinite-horizon optimal control problems, which are pieced together and exhibit a transient behavior. The turnpike property appears quite clearly in \cite{WildeKokotovic} for LQ problems under the Kalman condition, and is extended in \cite{AndersonKokotovic} to nonlinear control-affine systems where the vector fields are assumed to be globally Lipschitz, and being referred to as the exponential dichotomy property. In both cases the initial and final conditions for the trajectory are prescribed.
Their approach is remarkably simple and points out clearly the hyperbolicity phenomenon which is at the heart of the turnpike results. The use of Riccati-type reductions permits to quantify the saddle point property in a precise way.
Their proofs are however based on a Hamilton-Jacobi approach and, at the end of the article, the open question of extending their results to problems where the Hamilton-Jacobi theory cannot be used (that is, most of the time!) is formulated. In \cite{TrelatZuazua_JDE2015} we solved this open question by employing the Pontryagin maximum principle that yields a two-point boundary value problem.
In turn, in addition to the existing turnpike results at that time, the result of \cite{TrelatZuazua_JDE2015} gave a turnpike inequality that is also valid for the costate, having interesting consequences for numerical implementation (see Section \ref{sec_shooting} further).

\medskip
Recently there have been plenty of contributions on the turnpike, which appears to be an ubiquitous phenomenon that one can now identify in many fields. 
In what follows, we attempt to give an overview of the veritable explosion of novelties in this field, by describing some of the many generalizations and applications of which we are aware. Since, despite of many progresses on the topic, some simple examples still do escape to the analysis, we will also devote some time and room to analyze a meaningful nonlinear example in which global turnpike phenomena raise challenging questions.

\medskip

This article is structured as follows.

In Section \ref{sec_LQ}, we establish the exponential turnpike property for linear-quadratic (LQ) optimal control problems, in a very elementary way. The main result is Theorem \ref{thm_turnpike_LQ}.
Section \ref{sec_nonlinear} is devoted to showing how Theorem \ref{thm_turnpike_LQ} can be extended to finite-dimensional nonlinear optimal control problems, obtaining a local exponential turnpike property (i.e., the main result of \cite{TrelatZuazua_JDE2015}, improved with results of \cite{Trelat_MCSS2023}).
In Section \ref{sec_global}, we illustrate in detail the local and global aspects on simple but meaningful examples.
In Section \ref{sec_generalizations}, we provide a brief overview of some generalizations in various contexts, and explain why and how the turnpike phenomenon can be useful in some applications.

%\section{Exponential turnpike property for linear-quadratic optimal control problems}\label{sec_LQ}
\section{Linear-quadratic optimal control problems}\label{sec_LQ}
\subsection{Context and preliminary considerations}
Let $n$ and $m$ be nonzero integers. Let $A$, $B$, $Q$ and $U$ be real-valued matrices, of size, respectively, $n\times n$, $n\times m$, $n\times n$ and $m\times m$, with $Q$ and $U$ being symmetric and positive definite. 
Let $x_0, x_1, x_d\in\R^n$ and $u_d\in\R^m$ be fixed.
Given any $T>0$, we consider the linear-quadratic (LQ) optimal control problem
\begin{equation}\label{pbLQ}
\begin{split}
& \dot x(t) = Ax(t) + Bu(t), \\
& x(0)=x_0,\quad x(T)=x_1, \\
& \min \int_0^T \left[ (x(t)-x_d)^\top Q(x(t)-x_d) + (u(t)-u_d)^\top U (u(t)-u_d) \right] dt,
\end{split}
\end{equation}
where the state is $x(t)\in \R^n$ and the control is $u(t)\in \R^m$ (no constraint).

Assuming that $(A,B)$ satisfies the Kalman condition
\begin{equation}\label{rank}
\hbox{rank}\, (B, AB,...A^{n-1}B)=n,
\end{equation} 
i.e., that the above linear autonomous control system is controllable, by strict convexity there exists a unique optimal solution $(x_{_T}(\cdot),u_{_T}(\cdot))$ of \eqref{pbLQ} (see \cite{Kwakernaak, LeeMarkus, Trelat_book2005, Trelat_SB}).

Note that the controllability restriction is important here to assure that the set of admissible controls is non-void. Of course, if either the terminal or the initial condition is removed, the existence and uniqueness of the optimal control guaranteed for any pair $(A, B)$ without any further restriction. Yet, the controllability of $(A, B)$ plays an important role when aiming to guarantee the turnpike property.

It is well known that, by the Pontryagin maximum principle (see \cite{Kwakernaak, LeeMarkus, Pontryagin, Trelat_book2005, Trelat_SB, Vinter}), as a necessary condition for optimality, there must exist a map $\lambda_{_T}(\cdot):[0,T]\rightarrow\R^n$ called \emph{adjoint vector} or \emph{costate} such that 
$$
u_{_T}(t)=u_d+U^{-1}B^\top \lambda_{_T}(t)
$$
and
\begin{equation}\label{extr_syst1}
\begin{split}
\dot x_{_T}(t) &= A x_{_T}(t) + BU^{-1}B^\top \lambda_{_T}(t) + Bu_d, \\
\dot \lambda_{_T}(t) &= Qx_{_T}(t) - A^\top \lambda_{_T}(t) - Qx_d,
\end{split}
\end{equation}
for every $t\in[0,T]$,
i.e.,
\begin{equation}\label{extr_syst2}
\frac{d}{dt} \begin{pmatrix} x_{_T}(t)\\ \lambda_{_T}(t)\end{pmatrix}
= M \begin{pmatrix} x_{_T}(t)\\ \lambda_{_T}(t)\end{pmatrix} + \begin{pmatrix} Bu_d\\ -Qx_d\end{pmatrix},
\end{equation}
where
\begin{equation}\label{def_matrixM_LQ}
M = \begin{pmatrix}
A & BU^{-1}B^\top \\
Q & -A^\top
\end{pmatrix}.
\end{equation}
In this LQ setting, the adjoint state $\lambda_{_T}(\cdot)$ is unique and the above optimality system, together with the initial and terminal conditions on the state, namely, 
\begin{equation}\label{initial+final}
x(0)=x_0,\quad x(T)=x_1,
\end{equation} is a necessary and sufficient condition for optimality. 

We note that the system \eqref{extr_syst1}, together with the boundary conditions \eqref{initial+final}, does not fit into the classical framework of Cauchy problems for systems of differential equations, as no initial or terminal condition is imposed on the adjoint state $\lambda$, while both are specified for the state variable $x$. 
Although the existence and uniqueness of solutions to \eqref{extr_syst1}--\eqref{initial+final} follow directly from the fact that this system characterizes the optimal solution of the above optimal control problem, such well-posedness is not immediately evident from the perspective of dynamical systems theory.

\medskip
Under the above assumptions, the matrix $M$ is invertible (see further). Let $(\bar x,\bar\lambda)$ be the unique equilibrium point of the extremal system \eqref{extr_syst1} (equivalently, of \eqref{extr_syst2}), that is, the unique solution of 
\begin{equation}\label{extrLQ}
\begin{split}
A\bar x + BU^{-1}B^\top\bar\lambda + Bu_d &= 0 \\
Q\bar x - A^\top\bar\lambda - Q x_d &= 0,
\end{split}
\end{equation}
i.e., of the linear system
\begin{equation}\label{extrLQ_matrix}
M \begin{pmatrix} \bar x\\ \bar\lambda\end{pmatrix} + \begin{pmatrix} Bu_d\\ -Qx_d\end{pmatrix} = \begin{pmatrix} 0\\ 0\end{pmatrix},
\end{equation}
and set 
\begin{equation}\label{extrLQ2}
\bar u = u_d + U^{-1}B^\top\bar\lambda .
\end{equation}
The triple $(\bar x,\bar\lambda,\bar u)$ is the unique solution of the constrained optimization problem
\begin{equation}\label{staticpb_LQ} 
%\min_{(x,u)\in\R^n\times\R^m,\atop Ax+Bu=0} 
\min_{{(x,u)\in\R^n\times\R^m}\atop{Ax+Bu=0}}
\left[ (x-x_d)^\top Q(x-x_d) + (u-u_d)^\top U(u-u_d) \right],
\end{equation} 
and $\bar\lambda$ is the unique Lagrange multiplier (uniqueness is proved below).

\subsection{Main result: exponential turnpike property}
The \emph{exponential turnpike phenomenon}, that we formulate hereafter, stipulates that the triple $(x_{_T}(\cdot),\lambda_{_T}(\cdot),u_{_T}(\cdot))$, solution of the \emph{dynamic} optimization (optimal control) problem \eqref{pbLQ} (equivalently, of \eqref{extr_syst1}), remains exponentially close to the triple $(\bar x,\bar\lambda,\bar u)$, solution of the \emph{static} optimization problem \eqref{staticpb_LQ}, except at the beginning and at the end of the time interval $[0,T]$, something unavoidable in view of the initial and final conditions \eqref{initial+final}.

\begin{theorem}\label{thm_turnpike_LQ}
There exist constants $T_0>0$, $C>0$ and $\nu>0$, depending on $A,B,Q,U$ but not on $T,x_0,x_1$, such that, for every $T\geq T_0$, for all $x_0,x_1\in\R^n$, for every $t\in[0,T]$,
\begin{equation}\label{expturnLQ}
\Vert x_{_T}(t)-\bar x\Vert + \Vert u_{_T}(t)-\bar u\Vert + \Vert \lambda_{_T}(t)-\bar\lambda\Vert  
\leq C \Big( \Vert x_0-\bar x\Vert + \Vert x_1-\bar x\Vert \Big) \Big( e^{-\nu t} + e^{-\nu (T-t)} \Big)  .
\end{equation}
\end{theorem}

The exponential turnpike inequality \eqref{expturnLQ}, which implies that, except near $t=0$ and $t=T$, the triple (optimal state, optimal control, costate) at time $t$ is exponentially close to the steady-state solution of the static optimization problem \eqref{staticpb_LQ}, is of global nature. It holds for all initial and terminal constraints \eqref{initial+final} imposed on the optimal control problem, for all sufficiently large time horizons.

\begin{figure}[h]
%\begin{center}
%\resizebox{10cm}{!}{\input fig_turn.pdf_t}
%\end{center}
\centerline{\includegraphics[width=10cm]{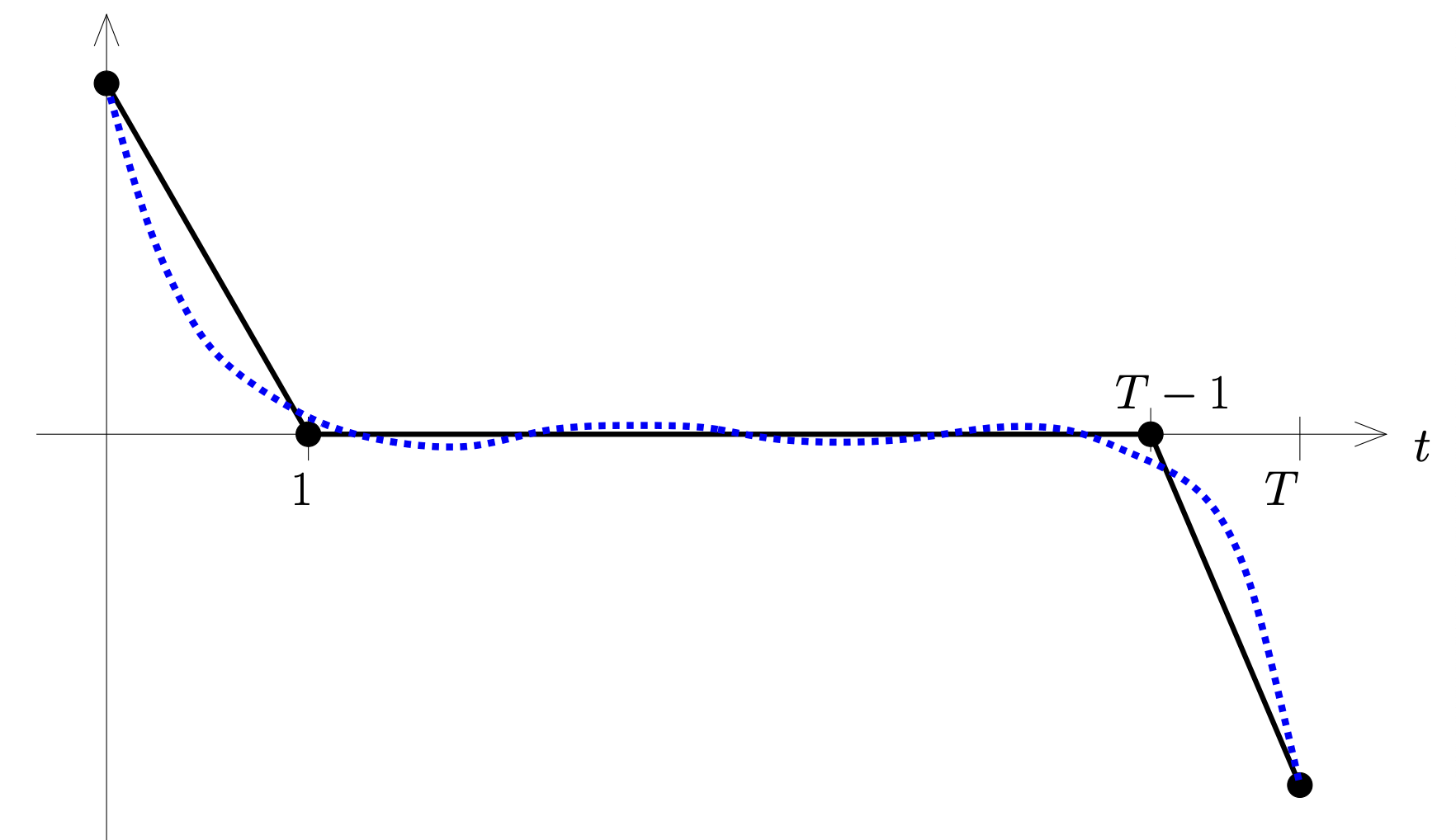}}
\caption{Illustration of the turnpike phenomenon: optimal strategy in dashed blue; quasi-optimal strategy in plain.}\label{fig_1}
\end{figure}

\medskip
It follows from this result that, when $T$ is large, a \emph{quasi-optimal control strategy} for the evolution system consists in:
\begin{enumerate}
\item steering the control system from the initial state $x_0$ to the (so-called) turnpike point $\bar x$, for instance in time $1$ (this is possible by the Kalman condition);
\item remaining at the steady-state $\bar x$, with the constant control $\bar u$, on the long time interval $[1,T-1]$;
\item steering the control system from the turnpike point $\bar x$ to the final state $x_1$ in time $1$ (this is possible by the Kalman condition).
\end{enumerate}
See Figure \ref{fig_1} for an illustration.
This three-step control strategy is quasi-optimal and the loss of optimality can be measured thanks to the discrepancy estimate \eqref{expturnLQ}. The first and the third arcs can be seen as \emph{transient arcs}. 

\subsection{Proof of Theorem \ref{thm_turnpike_LQ}}\label{sec_proof_thm_turnpike_LQ}
The proof is  elementary in the present LQ setting.
We proceed in several steps.

\medskip
\noindent {\bf Step 1. The static problem.} We start by proving the above claim that the static problem \eqref{staticpb_LQ} has a unique solution $(\bar x,\bar u)$.
Indeed, the Lagrangian of this optimization problem is 
$$
L(x,u,\lambda,\lambda^0) = \langle\lambda,Ax+Bu\rangle+\lambda^0\left( (x-x_d)^\top Q(x-x_d) + (u-u_d)^\top U(u-u_d) \right) .
$$
The first-order necessary conditions (which are also sufficient here by strict convexity) are obtained by differentiating $L$ with respect to $x$ and $u$. 
The qualification condition is satisfied because the map $(x,u)\mapsto Ax+Bu$ is surjective: indeed one has 
$$
\ker(A^\top)\cap\ker(B^\top)=\{0\}
$$ 
as a consequence of the Kalman condition \eqref{rank} (in fact the Kalman rank condition is equivalent to the fact that $\ker(B^\top)$ does not contain any eigenvector of $A^\top$). 
Therefore, we have $\bar\lambda^0 \neq 0$ (i.e., the abnormal case does not occur). Since the pair $(\bar\lambda, \bar\lambda^0)$ is defined only up to scaling and is unique up to a multiplicative constant, we can normalize it by setting $\bar\lambda^0 = -1/2$. With this normalization, the Lagrange multiplier $\bar\lambda$ is uniquely determined. 

The pair $(\lambda_{_T}(\cdot),\lambda_{_T}^0)\neq(0,0)$, defined up to scaling, has an interpretation in terms of Lagrange multipliers: indeed, the proof of the Pontryagin maximum principle relies on a preliminary application of the Lagrange multiplier rule to the optimal control problem viewed as a constrained (infinite-dimensional) optimization problem (see \cite{Trelat_book2005, TrelatJOTA} for a detailed explanation).

Recall that, in the time-dependent setting (with a natural static counterpart), an extremal is by definition a quadruple $(x_{_T}(\cdot), \lambda_{_T}(\cdot), \lambda^0_{_T}, u_{_T}(\cdot))$ (the so-called extremal lift of the optimal solution $(x_{_T}(\cdot),u_{_T}(\cdot))$) solution of the  first-order optimality system. It is said to be normal if $\lambda^0_{_T}<0$ and abnormal if $\lambda^0_{_T}=0$. Since the pair $(\lambda_{_T}(\cdot), \lambda^0_{_T})$ is defined up to scaling (an optimal solution may have several independent extremal lifts), a normal extremal can always be normalized so that $\lambda^0_{_T}=-1/2$.

\medskip

The conditions \eqref{extrLQ} and \eqref{extrLQ2} then follow from the Lagrange multiplier rule, which, in this case, provides both necessary and sufficient conditions for optimality.

\medskip
\noindent {\bf Step 2. Pontryagin maximum principle.} Before proceeding, we briefly outline the application of the Pontryagin Maximum Principle to problem~\eqref{pbLQ}. The Hamiltonian of the problem is
$$
H(x,\lambda,\lambda^0,u) = \langle\lambda,Ax\rangle+\langle\lambda,Bu\rangle + \lambda^0 \left[ (x-x_d)^\top Q(x-x_d) + (u-u_d)^\top U(u-u_d) \right] .
$$
Since $(x_{_T}(\cdot),u_{_T}(\cdot))$ is optimal, there must exist a costate $\lambda_{_T}(\cdot):[0,T]\rightarrow\R^n$ and a real number $\lambda_{_T}^0\leq 0$ such that $(\lambda_{_T}(t),\lambda_{_T}^0)\neq 0$ for any $t\in[0,T]$ and
\begin{equation*}
\begin{split}
\dot x_{_T}(t) &= \frac{\partial H}{\partial\lambda}(x_{_T}(t), \lambda_{_T}(t), \lambda^0_{_T}, u_{_T}(t)) = Ax_{_T}(t) + Bu_{_T}(t)  , \\
\dot\lambda_{_T}(t) &= -\frac{\partial H}{\partial x}(x_{_T}(t), \lambda_{_T}(t), \lambda^0_{_T}, u_{_T}(t)) = - A^\top \lambda_{_T}(t) - 2\lambda_{_T}^0Q(x_{_T}(t)-x_d) , \\
0 &= \frac{\partial H}{\partial u}(x_{_T}(t), \lambda_{_T}(t), \lambda^0_{_T}, u_{_T}(t)) = B^*\lambda_{_T}(t) + 2\lambda^0_{_T} U(u_{_T}(t)-u_d) .
\end{split}
\end{equation*}

Here, there is no abnormal extremal, thanks to the Kalman condition \eqref{rank} (see \cite{Pontryagin, Trelat_SB}). 
Indeed, otherwise, the condition $\partial H/\partial u=0$ yields $B^\top\lambda_{_T}(t)=0$, and by successive derivations and using the fact that $\dot \lambda_{_T}(t)=-A^\top \lambda_{_T}(t)$, we obtain $(A^kB)^\top\lambda_{_T}(t)=0$, for all $k \ge 0$ which raises a contradiction with the Kalman condition  since $\lambda_{_T}(t)\neq 0$. 
At this point we have used the classical equivalence between the Kalman rank condition, the controllability of the state equation and the observability of the adjoint system $\dot \lambda_{_T}(t) =  - A^\top \lambda_{_T}(t)$
ensuring that $\lambda \equiv 0$ as soon as $B^\top\lambda \equiv 0$.

Hence, we can normalize the adjoint vector so that $\lambda^0_{_T}=-1/2$. Now, the condition $\partial H/\partial u=0$ yields $u_{_T}(t)=u_d+U^{-1}B^\top \lambda_{_T}(t)$, and the resulting extremal system is \eqref{extr_syst1}, as claimed.
The uniqueness of $\lambda_{_T}(\cdot)$ follows from the absence of abnormal minimizers (see, e.g., \cite{CJT}).

\medskip
\noindent {\bf Step 3. Exponential turnpike.}   Let us now establish the exponential turnpike estimate \eqref{expturnLQ}.
It follows from \eqref{extr_syst2} and \eqref{extrLQ_matrix} that
$$
\frac{d}{dt} \begin{pmatrix} x_{_T}(t)-\bar x\\ \lambda_{_T}(t)-\bar\lambda\end{pmatrix}
= M \begin{pmatrix} x_{_T}(t)-\bar x\\ \lambda_{_T}(t)-\bar\lambda\end{pmatrix} 
$$
The key observation is that all eigenvalues of the matrix $M$ have a nonzero real part (actually, the number of unstable modes is equal to the number of stable modes, because $M$ is a Hamiltonian matrix).
Let us prove this fact. As mentioned above, as a consequence of the Kalman condition on $(A,B)$, we have
\begin{equation}\label{hautus}
\ker(A^\top- \xi I_n)\cap\ker(B^\top)=\{0\}\qquad\forall \xi\in\mathbb{C} .
\end{equation}
Let $(z_1,z_2)\in\mathbb{C}^n\times\mathbb{C}^n$ and $\mu\in\R$ such that $(M-i\mu)\begin{pmatrix}z_1\\z_2\end{pmatrix}=0$. Then,
\begin{align*}
(A-i\mu)z_1+BU^{-1}B^\top z_2&=0\\
Qz_1-(A^\top+i\mu) z_2&=0,
\end{align*}
hence $z_1=Q^{-1}(A^\top+i\mu)z_2$ and thus 
$$
(A-i\mu)Q^{-1}(A^\top+i\mu)z_2+BU^{-1}B^\top z_2=0 .
$$
Multiplying to the left by $\bar z_2^\top$, we obtain 
$$
\Vert Q^{-1/2}(A^\top+i\mu)z_2\Vert^2+\Vert U^{-1/2}B^\top z_2\Vert^2=0
$$
and hence $(A^\top+i\mu)z_2=0$ and $B^\top z_2=0$. We infer that $z_2=0$ by using \eqref{hautus}. We conclude that the matrix $M$ has no purely imaginary eigenvalues, as claimed. 

As a consequence, since $M$ is Hamiltonian, its stable and unstable subspaces $E^-$ and $E^+$ both have dimension $n$. We next use the boundary conditions on the state to estimate the corresponding stable and unstable components.

Let us first note that the projections of $E^-$ and $E^+$ onto the $x$-space are isomorphisms. Indeed, if $(0,\ell)^\top\in E^-$, the solution $(x(t),\lambda(t))^\top=e^{tM}(0,\ell)^\top$ tends to $0$ as $t\rightarrow+\infty$. Since
$$
\frac{d}{dt}\langle x(t),\lambda(t)\rangle = \Vert Q^{1/2}x(t)\Vert^2+\Vert U^{-1/2}B^\top\lambda(t)\Vert^2,
$$
and $x(0)=0$, integration on $[0,+\infty)$ gives $x\equiv 0$ and $B^\top\lambda\equiv 0$. By the Kalman condition, this implies $\ell=0$. The proof for $E^+$ is the same, integrating on $(-\infty,0]$. Thus there exist matrices $E_-$ and $E_+$ such that $E^-=\{(v,E_-v)^\top\,\mid\, v\in\R^n\}$ and $E^+=\{(w,E_+w)^\top\,\mid\, w\in\R^n\}$.

We can therefore write, in a unique way,
\begin{equation}\label{cdv}
\begin{pmatrix} x_{_T}(t)-\bar x\\ \lambda_{_T}(t)-\bar\lambda\end{pmatrix}
= \begin{pmatrix} v(t)+w(t)\\ E_-v(t)+E_+w(t)\end{pmatrix},
\end{equation}
where $(v(t),E_-v(t))^\top\in E^-$ and $(w(t),E_+w(t))^\top\in E^+$. By the exponential stability on $E^-$ and the exponential stability in reverse time on $E^+$, there exist constants $C_1,C_2,\nu_1,\nu_2>0$, depending only on $M$, such that
$$
\Vert v(t)\Vert \leq C_1 e^{-\nu_1 t}\Vert v(0)\Vert,
\qquad
\Vert w(t)\Vert \leq C_2 e^{-\nu_2 (T-t)}\Vert w(T)\Vert .
$$
Set $\eta_0=x_0-\bar x$ and $\eta_1=x_1-\bar x$. From the state boundary conditions we have
$$
v(0)+w(0)=\eta_0,
\qquad
v(T)+w(T)=\eta_1 .
$$
Using the preceding estimates at $t=T$ for $v$ and at $t=0$ for $w$, we get
$$
\Vert v(0)\Vert+\Vert w(T)\Vert \leq \Vert\eta_0\Vert+\Vert\eta_1\Vert + C_2e^{-\nu_2T}\Vert w(T)\Vert+C_1e^{-\nu_1T}\Vert v(0)\Vert .
$$
Choosing $T_0>0$ large enough so that $C_1e^{-\nu_1T}+C_2e^{-\nu_2T}\leq 1/2$ for every $T\geq T_0$, we obtain
$$
\Vert v(0)\Vert+\Vert w(T)\Vert \leq 2\big(\Vert x_0-\bar x\Vert+\Vert x_1-\bar x\Vert\big) .
$$
Combining the last estimate with \eqref{cdv} and with $u_{_T}-\bar u=U^{-1}B^\top(\lambda_{_T}-\bar\lambda)$ gives \eqref{expturnLQ}.

\begin{remark}\label{rem_ricc}
At the end of the above proof, one sees that the constants $T_0$, $C$ and $\nu$ of Theorem \ref{thm_turnpike_LQ} only depend on the hyperbolic splitting of $M$, and therefore only on the matrices, as stated.

The proof of Theorem \ref{thm_turnpike_LQ} that is given in \cite{TrelatZuazua_JDE2015} is a bit different and uses the algebraic Riccati theory. It yields an exponential turnpike inequality that is slightly more precise than \eqref{expturnLQ} (see \cite[Remark 13]{TrelatZuazua_JDE2015}). Actually, in this reference, the change-of-basis matrix $P$ is built by considering the minimal and maximal solutions of the Riccati algebraic equation, and the constants $C$ and $\nu$ are defined in terms of these matrices. 
\end{remark}

%\section{Generalizations}
%The turnpike phenomenon, that we have established in the previous section in an elementary way in the LQ case, has a number of generalizations that we describe hereafter. We first start with some historical comments.

%\section{Local exponential turnpike property for nonlinear optimal control problems}\label{sec_nonlinear}
\section{Finite-dimensional nonlinear optimal control problems}\label{sec_nonlinear}
The exponential turnpike given in Theorem \ref{thm_turnpike_LQ} in the LQ case was actually established in \cite{TrelatZuazua_JDE2015}, more generally, for finite-dimensional nonlinear optimal control problems, under appropriate assumptions.
The version that we derive hereafter is an improvement, borrowed from \cite{TrelatZuazua_ongoing}.

\subsection{Optimal control problem}
\paragraph{Setting.}
Let $n,m,k\in\N^*$ and let $f:\R^n\times\R^m\rightarrow\R^n$, $f^0:\R^n\times\R^m\rightarrow\R$ and $R %=(R_1,\ldots,R_k)
:\R^n\times\R^n\rightarrow\R^k$ be mappings of class $C^2$.
Let $\Omega$ be a measurable subset of $\R^m$.
Given any $T>0$, we consider the following nonlinear optimal control problem in $\R^n$, in fixed final time $T$: %, without control constraint:
\begin{align}
& \dot x(t) = f(x(t),u(t)) \label{syst} \\[1mm]
& R(x(0),x(T))=0, \label{terminalconditions} \\
&  u(t)\in\Omega, \label{syst_Omega} \\
%\qquad R(x(0),x(T))=0 \\
& \min \int_0^T f^0(x(t),u(t))\, dt, \label{mincost}
\end{align}
where $u\in L^\infty([0,T],\Omega)$ is the control.
The mapping $R$ stands for various possible boundary conditions (see Example \ref{ex_R} further).

We say that the pair $(x_{_T}(\cdot), u_{_T}(\cdot))$ is a \emph{globally optimal solution} of the optimal control problem \eqref{syst}--\eqref{terminalconditions}--\eqref{syst_Omega}--\eqref{mincost} if $u_{_T}(\cdot) \in L^\infty([0,T], \Omega)$, $x_{_T}(\cdot) : [0,T] \to \R^n$ is absolutely continuous, satisfies \eqref{syst} almost everywhere on $[0,T]$, and fulfills the boundary conditions \eqref{terminalconditions}. Additionally, the following inequality holds:
\[
\int_0^T f^0\big(x_{_T}(t), u_{_T}(t)\big)\, dt \leq \int_0^T f^0\big(x(t), u(t)\big)\, dt
\]
for any other admissible pair $(x(\cdot), u(\cdot))$ satisfying \eqref{syst}--\eqref{terminalconditions}--\eqref{syst_Omega}.

We refer to $(x_{_T}(\cdot), u_{_T}(\cdot))$ as a \emph{locally optimal solution} in $V \times U$, where $V \subset \R^n$ and $U \subset \R^m$ are open subsets containing $(x(\cdot), u(\cdot))$, if the above inequality is required to hold only for admissible trajectories $x(\cdot)$ taking values in $V$ and controls $u(\cdot)$ taking values in $U \cap \Omega$.

\medskip

We assume that there exists $T_0>0$ such that, for every $T\geq T_0$, there exists at least one (locally or globally) optimal solution $(x_{_T}(\cdot),u_{_T}(\cdot))$ of \eqref{syst}-\eqref{terminalconditions}-\eqref{syst_Omega}-\eqref{mincost}.

\medskip

Sufficient conditions ensuring existence are standard (see, e.g., \cite{Cesari, Trelat_SB}).
We do not assume uniqueness, but we note that uniqueness is ensured under differentiability properties of the value function and of the Hamiltonian (see, e.g., \cite[Theorem 7.3.9]{CannarsaSinestrari}), and thus is ``generic" in some sense). 

What we say hereafter is valid for locally or globally optimal solutions.

\paragraph{Application of the Pontryagin maximum principle.}
Denoting by $\langle\cdot,\cdot\rangle$ the Euclidean scalar product in $\R^n$ and by $\Vert\cdot\Vert$ the corresponding norm, we define the Hamiltonian of the optimal control problem \eqref{syst}-\eqref{terminalconditions}-\eqref{syst_Omega}-\eqref{mincost} by
\begin{equation}\label{def_Ham}
H(x,\lambda,\lambda^0,u) = \langle \lambda,f(x,u)\rangle + \lambda^0 f^0(x,u) = f(x,u)^\top\lambda + \lambda^0 f^0(x,u).
\end{equation}
By the Pontryagin maximum principle (see \cite{LeeMarkus, Pontryagin, Trelat_SB}), there exist $\lambda_{_T}^0\leq 0$ and an absolutely continuous mapping $\lambda_{_T}(\cdot):[0,T]\rightarrow\R^n$ (adjoint vector or costate), satisfying $(\lambda_{_T}(\cdot),\lambda_{_T}^0)\neq(0,0)$, such that
\begin{equation}\label{extremal_syst}
\begin{split}
& \dot x_{_T}(t) = \frac{\partial H}{\partial \lambda}(x_{_T}(t),\lambda_{_T}(t),\lambda_{_T}^0,u_{_T}(t)), \\
& \dot \lambda_{_T}(t) = -\frac{\partial H}{\partial x}(x_{_T}(t),\lambda_{_T}(t),\lambda_{_T}^0,u_{_T}(t)), \\
& u_{_T}(t) = \underset{v\in\Omega}{\mathrm{argmax}}\, H(x_{_T}(t),\lambda_{_T}(t),\lambda_{_T}^0,v) 
\end{split}
\end{equation}
for almost every $t\in[0,T]$.
We have moreover the transversality conditions on the adjoint vector at initial and final time: there exists $\gamma_{_T}\in\R^k$ such that
\begin{equation}\label{transversality}
\begin{pmatrix} -\lambda_{_T}(0)\\ \lambda_{_T}(T)\end{pmatrix} = dR(x_{_T}(0),x_{_T}(T))^\top \gamma_{_T} .
%\lambda_{_T}(0)\perp T_{x_{_T}(0)}M_0,\qquad \lambda_{_T}(T)\perp T_{x_{_T}(T)}M_1,
\end{equation}
%provided that the tangent space $T_{x_{_T}(0)}M_0$ to $M_0$ at $x_{_T}(0)$ and the tangent space $T_{x_{_T}(T)}M_1$ to $M_1$ at $x_{_T}(T)$ exist (these conditions are empty when the initial and final points are fixed in the optimal control problem).

\begin{example}\label{ex_R}
Let us give typical examples of mapping $R$ and of the corresponding transversality condition:
\begin{itemize}
\item For a fixed initial state $x_0$ and a fixed final state $x_1$, take $k=2n$ and $R(x,x')=(x-x_0,x-x_1)$. In this case, the transversality condition \eqref{transversality} gives no further information. 
\item For a fixed initial state $x_0$ and a free final state $x(T)$, take $k=n$ and $R(x,x')=x-x_0$. In this case the transversality condition \eqref{transversality} yields $\lambda_{_T}(T)=0$.
\item For a fixed initial state $x_0$ and a final state subject to the constraint $g(x(T))=0$ for some mapping $g=(g_1,\ldots,g_p):\R^n\rightarrow\R^p$ of class $C^2$ and some $p\in\N^*$, take $k=n+p$ and $R(x,x')=(x-x_0,g(x'))$. In this case the transversality condition \eqref{transversality}, with $\gamma_{_T} = (\hat\gamma_{_T}, \tilde\gamma_{_T})\in\R^n\times\R^k$, implies that $\lambda_{_T}(T)=dg(x_{_T}(T))^\top\tilde\gamma_{_T}$. % for some $\tilde\gamma_{_T}\in\R^k$.
\item For periodic conditions $x(0)=x(T)$, take $k=n$ and $R(x,x')=x-x'$. In this case the transversality condition \eqref{transversality} yields $\lambda_{_T}(0)=\lambda_{_T}(T)$. One can also consider partially periodic conditions.
\end{itemize}
\end{example}

\medskip

As mentioned above,  the pair $(\lambda_{_T}(\cdot),\lambda_{_T}^0)\neq(0,0)$ is defined up to scaling and,  
in view of \eqref{transversality}, it has an interpretation in terms of Lagrange multipliers. 
\medskip

The main result below (Theorem~\ref{thm_TZ}) establishes that, under suitable conditions, the optimal control problem admits a unique locally optimal solution $(x_{_T}(\cdot), u_{_T}(\cdot))$, which possesses a unique (up to scaling) normal extremal lift and satisfies the exponential turnpike property.

The analysis that we perform to show the turnpike property is valid for normal extremals but does not cover the case where an optimal solution $(x_{_T}(\cdot),u_{_T}(\cdot))$ would have only abnormal extremal lifts (but not any normal one). This abnormal situation is excluded by our assumptions. 

We recall that the existence of a normal extremal lift is ensured if there is no abnormal minimizer. This is so if, for instance, the initial state is fixed and the final state is let free (or the converse), as a consequence of \eqref{transversality} and of the nontriviality of the pair $(\lambda_{_T}(\cdot),\lambda^0)$. The absence of abnormal minimizers can be guaranteed under strong Lie bracket assumptions (see, e.g., \cite{AgrachevSachkov}) or for generic systems (see \cite{CJT}). 
It can also be noticed that the differentiability property of the value function (i.e., the minimal cost function between two points $x_0$ and $x_1$, in time $T$) ensures not only the uniqueness of the optimal trajectory between $x_0$ and $x_1$ in time $T$ but also the uniqueness of the extremal lift, which is moreover normal (see \cite[Theorem 7.3.9]{CannarsaSinestrari}).
When none of the known sufficient conditions apply, the normality assumption can be verified numerically by computing both normal and abnormal extremals and checking that, for the given boundary conditions, no abnormal extremal is minimizing (see~\cite{TrelatJOTA}). However, in some cases, abnormal minimizers do exist.

\begin{remark}\label{rem_int}
The last condition in \eqref{extremal_syst} gives
$$
\frac{\partial H}{\partial u}(x_{_T}(t),\lambda_{_T}(t),\lambda^0,u_{_T}(t)) = 0
$$
at every Lebesgue time $t\in[0,T]$ at which $u_{_T}(\cdot)$ is approximately continuous and $u_{_T}(t)\in\mathring{\Omega}$ (interior of $\Omega$).
\end{remark}

\subsection{Static optimization problem}\label{sec_static_nonlinear}
We consider the \textit{static} optimization problem (not depending on $T$)
\begin{equation}\label{staticpb}
\min \left\{ f^0(x,u) \ \mid\ (x,u)\in \R^n\times\Omega,\ f(x,u)=0 \right\}
\end{equation}
i.e., in other words, the problem of minimizing the instantaneous cost $f^0(x,u)$ over all possible equilibrium points of the dynamics.

\medskip

We assume that there exists at least one (local or global) minimizer $(\bar x,\bar u)\in\R^n\times\mathring{\Omega}$, i.e., with $\bar u$ in the interior of $\Omega$.

\medskip

Conditions ensuring the existence of a minimizer are standard in classical optimization (see, e.g., \cite{Bertsekas}).

As in the case of the optimal control problem, the results presented hereafter hold regardless of whether the minimizer $(\bar{x}, \bar{u})$ is local or global.

\paragraph{Application of the Lagrange multiplier rule.}Since $\bar u \in \mathring{\Omega}$, the Lagrange multiplier rule (see~\cite{Bertsekas}) ensures the existence of $\bar\lambda \in \R^n$ and $\bar\lambda^0 \leq 0$, with $(\bar\lambda, \bar\lambda^0) \neq (0,0)$, such that\footnote{We deliberately deviate from the standard convention in classical optimization, where $\lambda^0 \geq 0$ (as in~\cite{Bertsekas}), in order to express the optimality conditions in Hamiltonian form \eqref{staticextr}, using the Hamiltonian $H$ defined in \eqref{def_Ham}.}
\begin{equation*}%\label{LMstatic}
f(\bar x, \bar u) = 0 \qquad \text{and} \qquad
\bar\lambda^0 \, df^0(\bar x, \bar u)^\top + df(\bar x, \bar u)^\top \bar\lambda = 0
\end{equation*}
(the Lagrange multiplier condition). Equivalently, using the Hamiltonian $H$ defined in \eqref{def_Ham}, these conditions can be written as
\begin{equation}\label{staticextr}
\frac{\partial H}{\partial \lambda}(\bar x, \bar\lambda, \bar\lambda^0, \bar u) = 0, \qquad
\frac{\partial H}{\partial x}(\bar x, \bar\lambda, \bar\lambda^0, \bar u) = 0, \qquad
\frac{\partial H}{\partial u}(\bar x, \bar\lambda, \bar\lambda^0, \bar u) = 0.
\end{equation}

Note that the pair $(\bar\lambda,\bar\lambda^0)$ is defined up to scaling. The real number $\bar\lambda^0$, called abnormal multiplier, is the Lagrange multiplier associated with the function to be minimized. The vector $\bar\lambda$ is the Lagrange multiplier associated with the equality constraint.

\medskip
In classical constrained optimization, there are well known conditions (called qualification conditions) ensuring that we are in the normal case, such as Mangasarian-Fromovitz, Slater, etc (see, e.g., \cite{Bertsekas, BonnansShapiro}). Such conditions are generic. The simplest qualification condition is the assumption that $df(\bar x,\bar u)$ be surjective: this assumption, which will be done in Theorem \ref{thm_TZ} further, ensures that $(\bar x,\bar u)$ has no abnormal extremal lift.
But sometimes abnormal minimizers may happen: this is the case when $df(\bar x,\bar u)^\top$ has a nontrivial kernel.

In Theorem \ref{thm_TZ}, we will also make some assumptions implying that $(\bar x,\bar u)$ is a \emph{strict local minimizer} of the static optimization problem \eqref{staticpb}.

\subsection{Main result: local exponential turnpike property}
The turnpike phenomenon asserts that, when the final time $T$ is large, the optimal solution $(x_{_T}(\cdot), \lambda_{_T}(\cdot), u_{_T}(\cdot))$ remains, for most of the time interval $[0, T]$ -- except near the initial time $t = 0$ and the final time $t = T$ -- \emph{essentially close} to the optimal static point $(\bar x, \bar \lambda, \bar u) \in \R^n \times \R^n \times \Omega$. That is, roughly speaking,
\begin{equation}\label{roughturnpike}
x_{_T}(t) \simeq \bar x, \quad \lambda_{_T}(t) \simeq \bar \lambda, \quad u_{_T}(t) \simeq \bar u \qquad \forall t \in [\eta, T - \eta]
\end{equation}
for some $\eta > 0$.

A key observation, in light of Remark~\ref{rem_int}, is that the triple $(\bar x, \bar\lambda, \bar u)$ constitutes an equilibrium point of the extremal system \eqref{extremal_syst}. Accordingly, the informal property \eqref{roughturnpike} expresses that the (locally or globally) optimal trajectory $(x_{_T}(\cdot), \lambda_{_T}(\cdot), u_{_T}(\cdot))$, which solves the system \eqref{syst}--\eqref{terminalconditions}--\eqref{syst_Omega}--\eqref{mincost}--\eqref{extremal_syst}--\eqref{transversality}, remains essentially close to this \emph{equilibrium point} of the dynamics \eqref{extremal_syst}. Moreover, this equilibrium $(\bar x, \bar\lambda, \bar u)$ is itself a (locally or globally) optimal solution of the static problem \eqref{staticpb}. This is the essence of the turnpike phenomenon.

\paragraph{Heuristic explanation of the turnpike property.}We now provide a heuristic explanation for why such a property is expected. Since $T$ is assumed to be large, we introduce the small parameter $\varepsilon = 1/T$. Next, we perform a time reparametrization by setting $s = \varepsilon t$, so that as $t$ ranges over $[0, T]$, the new variable $s$ spans $[0, 1]$. 

Defining $\mathrm{x}_\varepsilon(s) = x_{_T}(t)$, $\mathrm{\lambda}_\varepsilon(s) = \lambda_{_T}(t)$, and $\mathrm{u}_\varepsilon(s) = u_{_T}(t)$, the extremal system can be formally rewritten as:
\begin{equation}\label{informal}
\begin{split}
& \frac{\partial H}{\partial \lambda}\big(\mathrm{x}_\varepsilon(s), \mathrm{\lambda}_\varepsilon(s), -1, \mathrm{u}_\varepsilon(s)\big) = \varepsilon \mathrm{x}_\varepsilon'(s) \simeq 0, \\
& \frac{\partial H}{\partial x}\big(\mathrm{x}_\varepsilon(s), \mathrm{\lambda}_\varepsilon(s), -1, \mathrm{u}_\varepsilon(s)\big) = -\varepsilon \mathrm{\lambda}_\varepsilon'(s) \simeq 0, \\
& \frac{\partial H}{\partial u}\big(\mathrm{x}_\varepsilon(s), \mathrm{\lambda}_\varepsilon(s), -1, \mathrm{u}_\varepsilon(s)\big) = 0,
\end{split}
\end{equation}
where the last condition holds assuming that the control lies in the interior of $\Omega$.

These formal computations provide intuition for the turnpike property \eqref{roughturnpike}, particularly when $(\bar x, \bar \lambda, \bar u)$ is the unique solution of \eqref{staticextr} and under suitable nondegeneracy assumptions on the system \eqref{staticextr}. 

However, this argument remains heuristic. In particular, justifying \eqref{informal} rigorously requires showing that the derivatives $\mathrm{x}_\varepsilon'(s)$ and $\mathrm{\lambda}_\varepsilon'(s)$ are uniformly bounded (away from neighborhoods of $t = 0$ and $t = T$) independently of $\varepsilon$. Equivalently, this amounts to proving that $T \dot{x}_{_T}(t)$ and $T \dot{\lambda}_{_T}(t)$ remain uniformly bounded (again, except near $t = 0$ and $t = T$) independently of $T$. Controlling these derivatives precisely is the main technical challenge in establishing the turnpike property rigorously.

\paragraph{Main result.}
We use hereafter the notations
$$
\bar H_{\star\#} = \frac{\partial^2 H}{\partial\star\partial\#}(\bar x,\bar \lambda,-1,\bar u)
$$
with $\star,\#$ being either $x$ or $u$ or $\lambda$ ($\bar\lambda^0=-1$ will be ensured by the assumptions done in the theorem).
We set
$$
A = \bar H_{\lambda x} = \frac{\partial f}{\partial x}(\bar x,\bar u),\qquad
\bar B = \bar H_{\lambda u} = \frac{\partial f}{\partial u}(\bar x,\bar u), \qquad \bar U=-\bar H_{uu}
$$
and, since we assume hereafter that $\bar U$ is invertible,
$$
\bar A = A + \bar B \bar U^{-1} \bar H_{ux}, \qquad \bar W = -\bar H_{xx} - \bar H_{xu}\bar U^{-1}\bar H_{ux} .
$$
The next theorem is the main result of \cite{TrelatZuazua_JDE2015}, improved with some results of \cite{Trelat_MCSS2023}.%\footnote{In \cite{TrelatZuazua_JDE2015}, Theorem \ref{thm_TZ} is established under the additional \emph{smallness condition}

\begin{theorem}\label{thm_TZ}
We assume that there exists at least one (local or global) minimizer $(\bar x,\bar u)\in\R^n\times\mathring{\Omega}$ (i.e., with $\bar u$ in the interior of $\Omega$) of the static optimization problem \eqref{staticpb}, with an extremal lift $(\bar x,\bar\lambda,\bar\lambda^0,\bar u)$. Using the above notations, we make the following assumptions: 
\begin{enumerate}[label=(\roman*)]
\item\label{H1_U_W} $\bar U$ and $\bar W$ are positive definite symmetric matrices.
%\item\label{H1_W} $\bar W$ is a positive definite symmetric matrix.
\item\label{H1_Kalman} The pair $(A,\bar B)$ satisfies the Kalman condition\footnote{Equivalently, the pair $(\bar A,\bar B)$ satisfies the Kalman condition (indeed, $(A,\bar B)$ satisfies the Kalman condition if and only if $(A+\bar BK\bar B,\bar B)$ satisfies the Kalman condition for every matrix $K$ of size $m\times n$). \\
This ensures that the static optimization problem \eqref{staticpb} is qualified and one can take $\bar\lambda^0=-1$.}, i.e., $\mathrm{rank}(\bar B,A\bar B,\ldots,A^{n-1}\bar B)=n$. In other words, the linearized system at $(\bar x,\bar u)$ is controllable.
\item\label{R_nonsing} The point $(\bar x,\bar x)$ is not a singular point of the mapping $R$. %; moreover, the norm of the Hessian of $R$ at $(\bar x,\bar x)$ or the mapping $R$ is generic (see \cite{TrelatZuazua_JDE2015}).
\end{enumerate}
Then, we have the following results:
\begin{itemize}
\item {\bf Static optimization problem:}

The pair $(\bar x,\bar u)$ is a strict local minimizer of the static optimization problem \eqref{staticpb}. Moreover, it has a unique (up to scaling) extremal lift $(\bar x,\bar\lambda,\bar\lambda^0,\bar u)$, which is moreover normal. Hence, we can normalize it so that $\bar\lambda^0=-1$.
\item {\bf Local exponential turnpike property:}

There exist an open neighborhood $V_{\bar x}$ of $\bar x$ in $\R^n$ and an open neighborhood $U_{\bar x}$ of $\bar u$ in $\R^m$, and there exist $C>0$, $\nu>0$ and $T_0>0$ such that,
%There exist $\varepsilon>0$, $C>0$, $\nu>0$ and $T_0>0$ such that if 
%\begin{equation}\label{smallness}
%\Vert R(\bar x,\bar x)\Vert + \left\Vert \begin{pmatrix} -\bar \lambda\\ \bar \lambda\end{pmatrix} - dR(\bar x,\bar x)^\top \gamma_{_T} \right\Vert + \Vert d^2R(\bar x,\bar x)\Vert\leq\varepsilon
%\end{equation}
%then, 
for every $T\geq T_0$, the optimal control problem \eqref{syst}-\eqref{terminalconditions}-\eqref{syst_Omega}-\eqref{mincost} has a unique locally optimal solution $(x_{_T}(\cdot),u_{_T}(\cdot))$
in $V_{\bar x}\times U_{\bar x}$,
which moreover has a unique normal extremal lift $(x_{_T}(\cdot),\lambda_{_T}(\cdot),-1,u_{_T}(\cdot))$ solution of \eqref{extremal_syst}-\eqref{transversality} and satisfying
\begin{equation}\label{turnpikeexp}
\Vert x_{_T}(t)-\bar x\Vert + \Vert \lambda_{_T}(t)-\bar \lambda\Vert + \Vert u_{_T}(t)-\bar u\Vert \leq C \big( e^{-\nu t} + e^{-\nu(T-t)}\big) \qquad \forall t\in [0,T] .
\end{equation}
\end{itemize}
\end{theorem}

The exponential turnpike estimate \eqref{turnpikeexp} implies that, except near $t=0$ and $t=T$, the locally optimal triple $(x_{_T}(\cdot),\lambda_{_T}(\cdot),u_{_T}(\cdot))$ solution of \eqref{syst}-\eqref{terminalconditions}-\eqref{syst_Omega}-\eqref{mincost}-\eqref{extremal_syst}-\eqref{transversality} %, whose existence is stated in the theorem under the important smallness condition \eqref{smallness}, 
is exponentially close to the (locally or globally) optimal triple $(\bar x,\bar \lambda,\bar u)$ solution of \eqref{staticpb}-\eqref{staticextr}.
The constants $C$ and $\nu$ (not depending on $T$) can be made more precise by considering Riccati algebraic equations, as said in Remark \ref{rem_ricc}.

\begin{remark}
Note that the triple $(\bar x,\bar u,\bar\lambda)$ is an equilibrium point of the first-order optimality system \eqref{extremal_syst} given by the Pontryagin maximum principle.
The proof of the theorem consists, first, of linearizing the optimality system \eqref{extremal_syst} at the equilibrium $(\bar x,\bar u,\bar\lambda)$, as it was done in \cite{TrelatZuazua_JDE2015}, leading to a local turnpike property under an appropriate smallness condition, namely, that
$\Vert R(\bar x,\bar x)\Vert + \Vert( -\bar \lambda \ \, \bar \lambda)^\top - dR(\bar x,\bar x)^\top \gamma_{_T} \Vert + \Vert d^2R(\bar x,\bar x)\Vert\leq\varepsilon $
for some small enough $\varepsilon>0$, which is stronger than what is stated in Theorem \ref{thm_TZ} where the conclusion is established in a neighborhood $V_{\bar x}\times U_{\bar u}$ of $(\bar x,\bar u)$, not requiring any smallness condition on the adjoint vector.
Actually, under Assumptions \ref{H1_U_W}-\ref{H1_Kalman}-\ref{R_nonsing}, the nonlinear optimal control problem \eqref{syst}-\eqref{terminalconditions}-\eqref{syst_Omega}-\eqref{mincost} can be approximated by a LQ optimal control problem, near the local minimizer $(\bar x,\bar u)$ of the static optimization problem \eqref{staticpb}. Then, following the arguments of the proof of Theorem \ref{thm_turnpike_LQ}, one obtains a \emph{local} turnpike result, under the above smallness condition. 
The extension done in \cite{TrelatZuazua_ongoing} allows to drop this smallness assumption, by using sensitivity analysis and conjugate point theory (as in \cite[Proposition 1]{Trelat_MCSS2023}).

Theorem \ref{thm_TZ} thus gives a \emph{local} exponential turnpike property, in a neighborhood $V_{\bar x}\times U_{\bar x}$ of the pair $(\bar x,\bar u)$ that is a local or global minimizer of the static problem \eqref{staticpb}. 
It is noticeable that in $V_{\bar x}\times U_{\bar x}$ we have a unique (locally) optimal solution of the optimal control problem in time $T$ large enough. 
\end{remark}

%\begin{remark}
%Note that, in addition to \cite{TrelatZuazua_JDE2015}, we have put a control constraint $u(t)\in\Omega$ but we have assumed that $\bar u\in\mathring{\Omega}$, so that, in the long middle part of the time interval, the optimal control, which is close to $\bar u$, does not saturate the constraints. Bang arcs can only occur at the beginning and at the end of the time frame.
%
%We can as well add some state constraints $x(t)\in C$, without many significant changes, provided we assume that $\bar x\in\mathring{C}$, and that we use a version of the Pontryagin maximum principle with state constraints (see \cite{Vinter}) to describe the beginning and the end of the optimal trajectory.
%
%When $\bar u\in\partial\Omega$ and/or $\bar x\in\partial C$, i.e., when the optimal steady-state saturates the contraints, we do not know how to derive an exponential turnpike property by exploiting the Pontryagin extremal system, because our approach (described next) relies on linearizing this first-order optimality system at the optimal steady-state and requires smoothness properties of the extremal control. In such a case, the approach by local dissipativity, which is softer, is more appropriate but leads to weaker turnpike results, like the ``measure-turnpike property" (see Section \ref{sec_dissipativity}). 
%\end{remark}

\begin{remark}\label{rem4}
In the context of Theorems \ref{thm_turnpike_LQ} and \ref{thm_TZ}, taking $(x_0,x_1)\neq(\bar x,\bar x)$, the turnpike is \emph{never} reached exactly, in the sense that the trajectory $(x_{_T}(\cdot),\lambda_{_T}(\cdot))$ never crosses the point $(\bar x,\bar\lambda)$ (it  follows from the proof, by Cauchy uniqueness). 
Actually, $(x_{_T}(\cdot),\lambda_{_T}(\cdot))$ evolves exponentially close to $(\bar x,\bar\lambda)$ (see the numerical simulations provided further on several examples).
Accordingly, under analyticity assumptions, the optimal trajectory $x_{_T}(\cdot)$ has no nontrivial arc along which one would have $x_{_T}(t)\equiv\bar x$. 

There are other cases where turnpike occurs, not obtained by linearization as in Theorem \ref{thm_TZ}, where an ``exact turnpike" may happen (see the generalizations in Section \ref{sec_generalizations}). This is the case when, for example, the turnpike consists of a singular arc, which is a part of the optimal strategy. 

An explicit example is given by the celebrated Dubins problem (see \cite{Dubins}), in which it is well known that the optimal strategy between two points $x_0$ and $x_1$ consists of three arcs: the first and the third are circular arcs, and the middle arc is a straight. Actually, the middle arc is a singular arc corresponding to the zero control. Now, taking $x_0$ and $x_1$ far from each other, this forces the middle arc to be in long time, thus showing an \emph{exact} turnpike phenomenon, where the turnpike consists of the singular arc. Of course, here, the turnpike is not restricted to a singleton but consists of a curve (see Section \ref{sec_moregeneral} further).

Another example is in orbit transfer where the classically used strategy consists of performing a first brief thrust arc, then a long-time ballistic arc (i.e., with no thrust) and, at the end, another brief thrust arc used to reach the target. This is also a turnpike picture.
\end{remark}

\subsection{Global considerations}\label{sec_global}
As explained previously, Theorem \ref{thm_TZ} provides a \emph{local} exponential turnpike property. %under an appropriate smallness condition \eqref{smallness}. 
We will see at the end of Section \ref{sec_dissipativity} that, under additional assumptions (namely, uniqueness of minimizers of the optimal control problem and of the static problem, and a global dissipativity assumption), Theorem \ref{thm_TZ} can be adapted to get a \emph{global} exponential turnpike property. But global dissipativity is a very strong assumption, difficult (not to say impossible) to check on practical examples. 

Since Theorem~\ref{thm_TZ} is a local result, it is natural that when the static optimization problem \eqref{staticpb} admits multiple local or global minimizers, corresponding local or global turnpike properties may arise. For instance, if the static problem has two global minimizers, the state space may be partitioned into regions: in one region, globally optimal solutions of the dynamic optimal control problem exhibit a turnpike behavior associated with the first minimizer, while in another region, the turnpike is determined by the second minimizer.

Such interesting competition phenomena have been studied in~\cite{Rapaport, Rapaport2}, as discussed in Section~\ref{sec_moregeneral}. Additional numerical simulations illustrating these effects are provided in~\cite[Appendix]{Trelat_MCSS2023}. In this section, we present several further numerical simulations that give rise to new open questions.

\subsubsection{A first example}\label{sec_ex1}
\begin{figure}[H]
\centerline{\includegraphics[width=6cm]{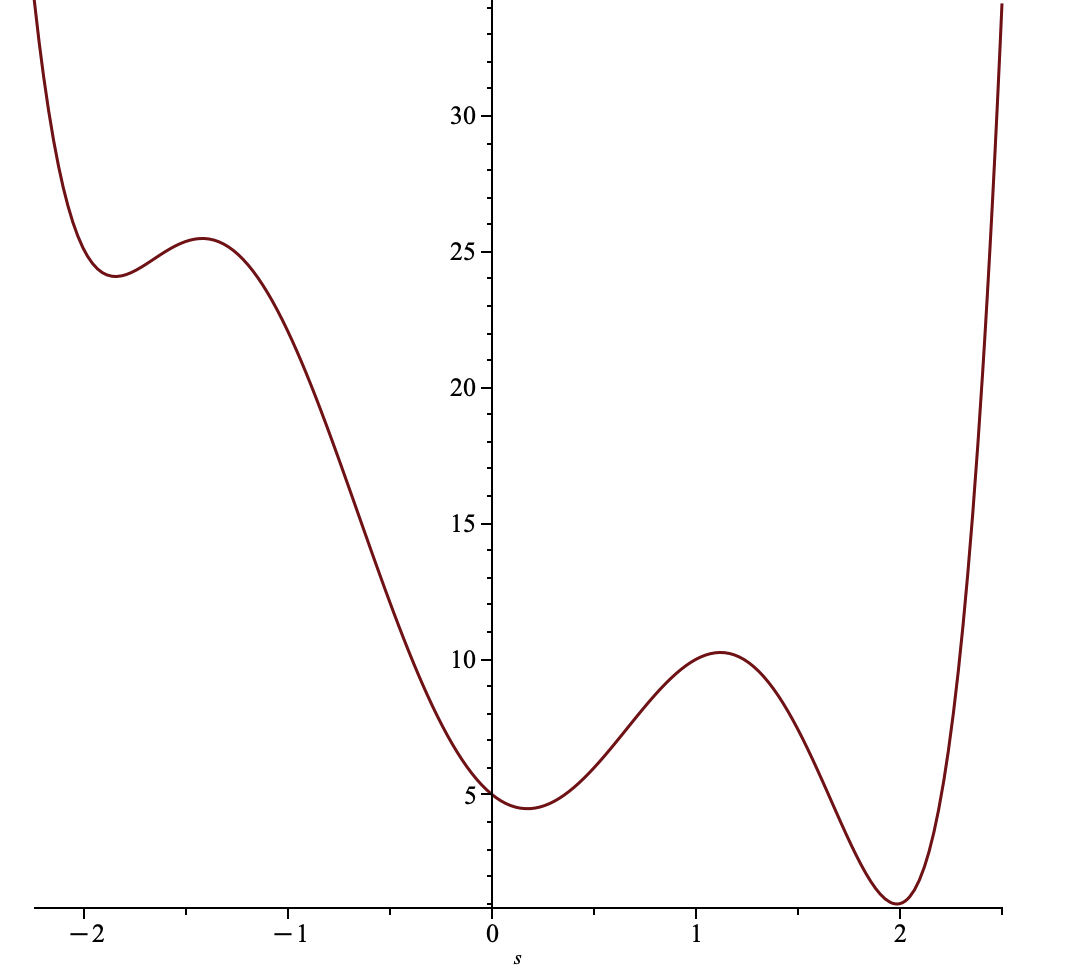}}
\caption{Graph of $s\mapsto(s-1)^2+(s-2)^2+(4s-s^3)^2$.}\label{fig_static_xcube}
\end{figure}

We consider the optimal control problem in $\R^2$: % (given in \cite{TrelatZuazua_JDE2015}):
\begin{equation}\label{exa}
\begin{split}
& \dot x(t) = x(t)-y(t), \\
& \dot y(t) = -4x(t)+y(t)^3 + u(t) , \\
& x(0)=x_0,\ y(0)=y_0, \\
& x(T)=x_1,\ y(T)=y_1,\\
& \min \int_0^T \big[ (x(t)-1)^2+(y(t)-2)^2+u(t)^2 \big]\, dt . 
\end{split}
\end{equation}
The static optimization problem
$$
\min_{\stackrel{x-y=0}{-4x+y^3+u=0}}  \big[ (x-1)^2+(y-2)^2+u^2 \big] 
$$
has a unique global minimizer given by
$$
\bar x^{\textrm{glob}}=\bar y^{\textrm{glob}}=1.98432,\quad\bar u^{\textrm{glob}}=0.123985,\quad\bar \lambda^{\textrm{glob}}=(2.96051,0.247969) , 
$$
and two local minimizers given by
$$
\bar x^{\textrm{loc}1}=\bar y^{\textrm{loc}1}=-1.84987,\quad\bar u^{\textrm{loc}1}=-1.06922,\quad\bar \lambda^{\textrm{loc}1}=(-14.2535,-2.13844) , 
$$
and
$$
\bar x^{\textrm{loc}2}=\bar y^{\textrm{loc}2}=0.171094,\quad\bar u^{\textrm{loc}2}=0.679368,\quad\bar \lambda^{\textrm{loc}2}=(3.77714,1.35874) .
$$
Solving the static problem is equivalent to minimizing the function $s\mapsto(s-1)^2+(s-2)^2+(4s-s^3)^2$, which is drawn on Figure \ref{fig_static_xcube}.

Assumptions \ref{H1_U_W} and \ref{H1_Kalman} of Theorem \ref{thm_TZ} are checked numerically for the three minimizers. Assumption \ref{R_nonsing} is obviously satisfied since $R((x,y),(x',y'))=(x-x_0,y-y_0,x'-x_1,y'-y_1)$.

\begin{figure}[H]
\centerline{\includegraphics[width=15cm]{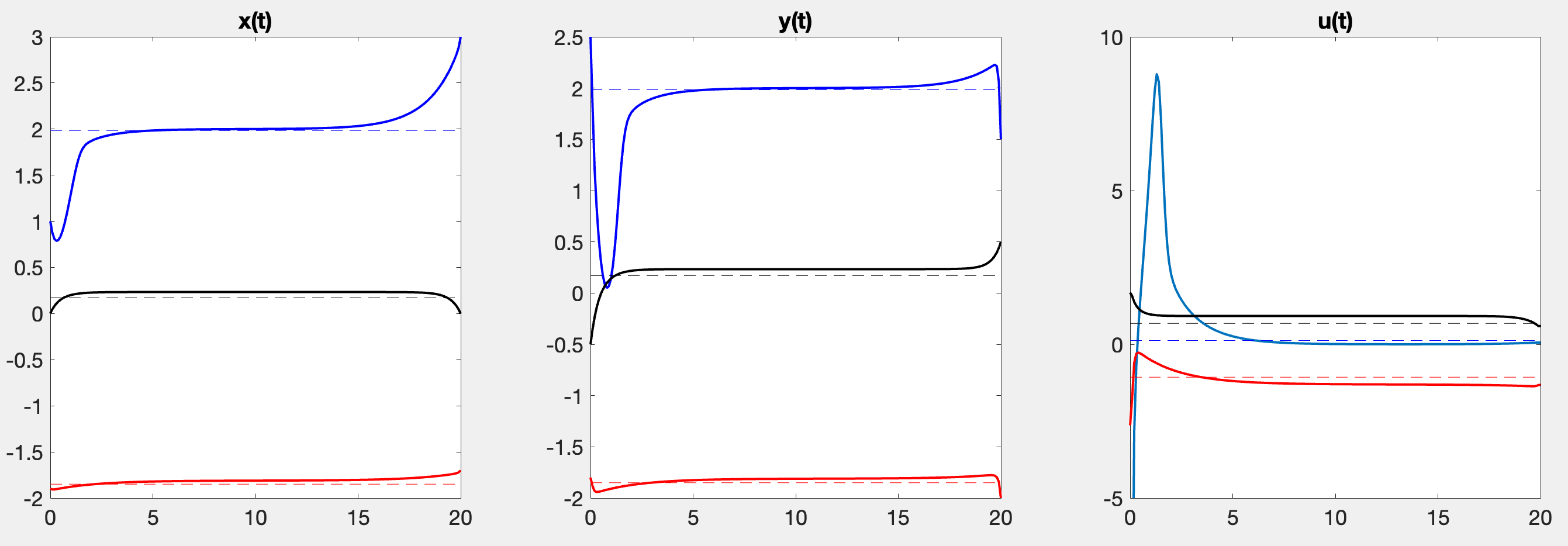}}
\caption{Locally optimal trajectories near the global static minimizer (in blue), and near the two local static minimizers (in red and black).}\label{fig_simu_xcube}
\end{figure}

\medskip
Figure \ref{fig_simu_xcube} gives some numerical simulations. 

The blue curve represents the locally optimal trajectory steering the control system in time $T = 20$ from $(x_0, y_0) = (1, 2.5)$ to $(x_1, y_1) = (3, 1.5)$. This trajectory exhibits the local exponential turnpike property around the global minimizer $(\bar x^{\textrm{glob}}, \bar y^{\textrm{glob}}, \bar u^{\textrm{glob}})$ of the static problem.

The red and black curves correspond, respectively, to locally optimal trajectories steering the system in time $T = 20$ from $(x_0, y_0) = (-1.9, -1.8)$ to $(x_1, y_1) = (-1.7, -2)$, and from $(x_0, y_0) = (0, 0)$ to $(x_1, y_1) = (-0.5, 0.5)$. Both trajectories display the local exponential turnpike property around their respective local minimizers $(\bar x^{\textrm{loc1}}, \bar y^{\textrm{loc1}}, \bar u^{\textrm{loc1}})$ and $(\bar x^{\textrm{loc2}}, \bar y^{\textrm{loc2}}, \bar u^{\textrm{loc2}})$ of the static problem.

Numerical computations confirm that the red and black trajectories are locally but not globally optimal. The globally optimal trajectories corresponding to these boundary conditions exhibit a global turnpike behavior around the global minimizer $(\bar x^{\textrm{glob}}, \bar y^{\textrm{glob}}, \bar u^{\textrm{glob}})$.

Similarly, numerical results confirm that the blue trajectory enjoys a global exponential turnpike property around the global minimizer $(\bar x^{\textrm{glob}}, \bar y^{\textrm{glob}}, \bar u^{\textrm{glob}})$.

This observation motivates the global perspective developed in the next section.

\subsubsection{Local or global turnpike}
At this step, the following natural questions emerge:
\begin{enumerate}
\item Assuming that the static optimization problem has a unique local minimizer, which is thus global, do the optimal trajectories enjoy a global exponential turnpike property around this minimizer? %(without invoking any global dissipativity assumption)
\item Assuming that the static optimization problem has a unique global minimizer and has at least another local minimizer, do the optimal trajectories enjoy a global exponential turnpike property around this global minimizer? %(without invoking any global dissipativity assumption)
\item Assuming that the static optimization problem has several global minimizers, i.e., assuming that we have several global turnpikes, is one of them ``better" than the others in some sense?
\end{enumerate}
It is proved in \cite{TrelatZuazua_ongoing} that the answers to the above three questions are all negative. 
Let us report on some examples and considerations done in this work.

\paragraph{Counterexample to global turnpike, in spite of a unique global minimizer of the static problem.}
In $\R^2$, we denote the canonical coordinates by $(x,y)$.
Given any $r>0$, we denote by $B(0,r)=\{(x,y)\in\R^2\ \mid\ x^2+y^2<r\}$ and by $\bar B(0,r)$ its closure.
Let $\varphi:[0,+\infty)\rightarrow\R$ be a function of class $C^\infty$, such that $0\leq\varphi\leq 1$,
\begin{equation*}
\varphi(s) = \left\{ \begin{array}{ll}
1 & \textrm{if}\ 0\leq s\leq 1, \\
0 & \textrm{if}\ s\geq 4 ,
\end{array}\right.
\end{equation*}
and $\varphi$ is decreasing on $[1,4]$.
Let $R>2$ to be chosen later. We set
$$
f^0(x,y,u) = \varphi(x^2+y^2) \left( (x-1)^2+y^2 \right) + (1-\varphi(x^2+y^2)) \left( x^2+y^2-R^2\right)^2 + u^2 .
$$
Let $x_0,x_1,y_0,y_1\in\R$. Given any $T>0$, we consider the optimal control problem in $\R^2$
\begin{equation}\label{syst_circle}
\begin{split}
&\dot x(t)= y(t), \qquad\qquad\qquad x(0)=x_0,\quad x(T)=x_1, \\
& \dot y(t)=-x(t)+u(t), \qquad\ y(0)=y_0,\quad y(T)=y_1,  \\
& \min \int_0^T f^0(x(t),y(t),u(t))\, dt, 
\end{split}
\end{equation}
with no constraint on the scalar control $u$.

It is proved in~\cite{TrelatZuazua_ongoing} that:
\begin{itemize}
    \item For $R$ sufficiently large, the static optimization problem associated with \eqref{syst_circle} admits a unique local (and thus global) minimizer, given by $\bar x = 1/2$, $\bar y = 0$, and $\bar u = 1/2$.
    
    \item Since the problem reduces to a linear-quadratic (LQ) optimal control problem inside the ball $B(0,1)$, it exhibits a turnpike property at $(\bar x, \bar y, \bar u) = \left(1/2, 0, 1/2\right)$. Consequently, the optimal control problem \eqref{syst_circle} satisfies the local exponential turnpike property near this triple.
    
    \item Outside the ball $B(0,2)$, the cost functional to be minimized is
    \[
    \int_0^T \left[ \left(x(t)^2 + y(t)^2 - R^2\right)^2 + u(t)^2 \right] dt,
    \]
    which vanishes along the trajectory
    \[
    x(t) = R \sin(t + \tau), \quad y(t) = R \cos(t + \tau), \quad u(t) = 0,
    \]
    for any $\tau \in \R$. This trajectory is a solution of the control system since $\dot{x}(t) = y(t)$ and $\dot{y}(t) = -x(t)$, and it traces the circle $\mathcal{C} = \{ (x, y) \in \R^2 \mid x^2 + y^2 = R^2 \}$, fully contained in the region $\R^2 \setminus B(0,2)$.
\end{itemize}
Therefore, there is no global turnpike property at $(\bar x,\bar y,\bar u)$: given any initial point $(x_0,y_0)$ and final point $(x_1,y_1)$ in $\R^2$, the globally optimal trajectory will always spend most of its time along $\mathcal{C}$ (except at the beginning and at the end).

\begin{figure}[H]
\centerline{\includegraphics[width=12cm]{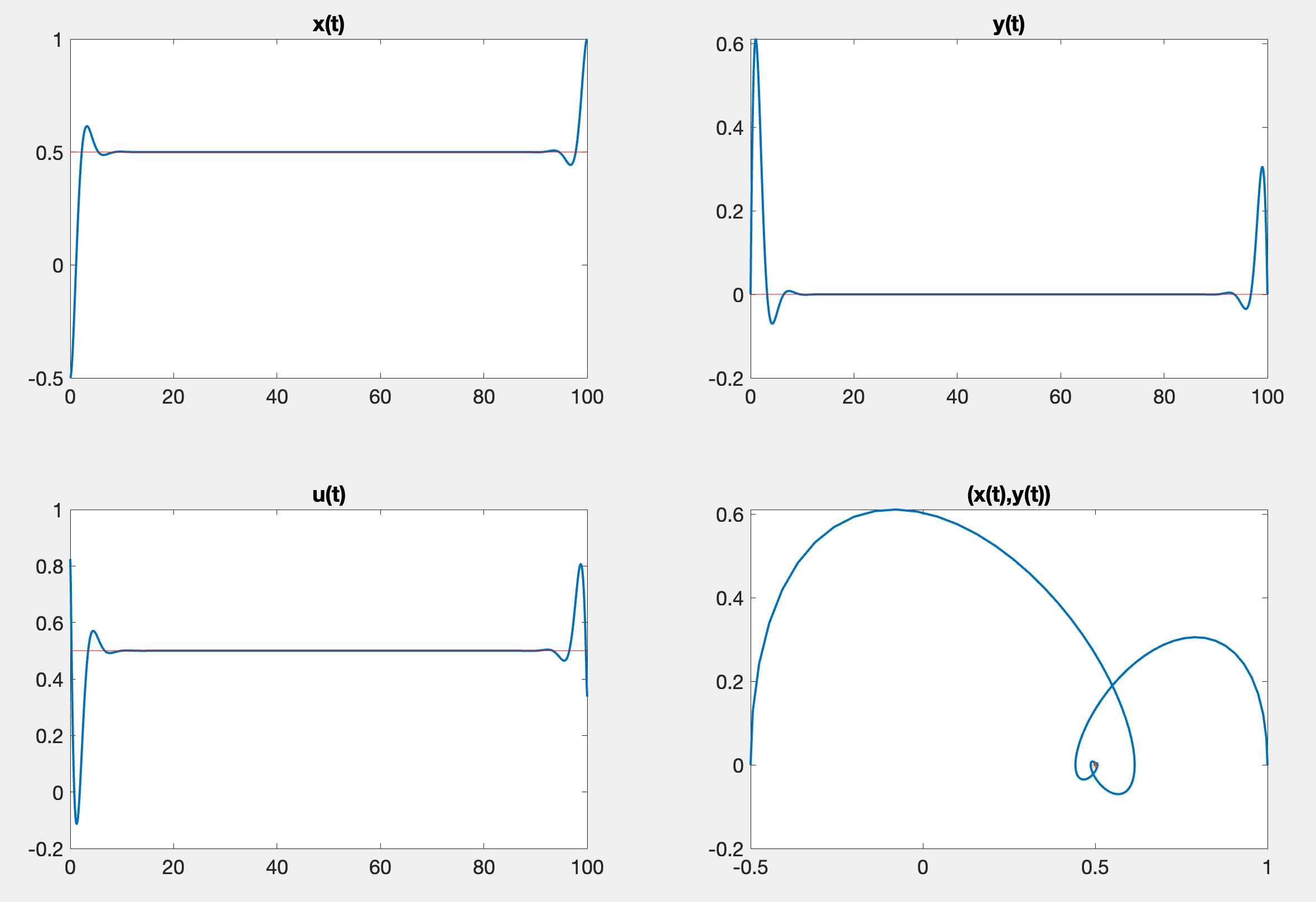}}
\caption{Locally optimal trajectory steering $(x_0,y_0)$ to $(x_1,y_1)$ in time $T$, of cost $\simeq 52.3983$.}\label{fig_simu_circle1}
\end{figure}

\begin{figure}[H]
\centerline{\includegraphics[width=12cm]{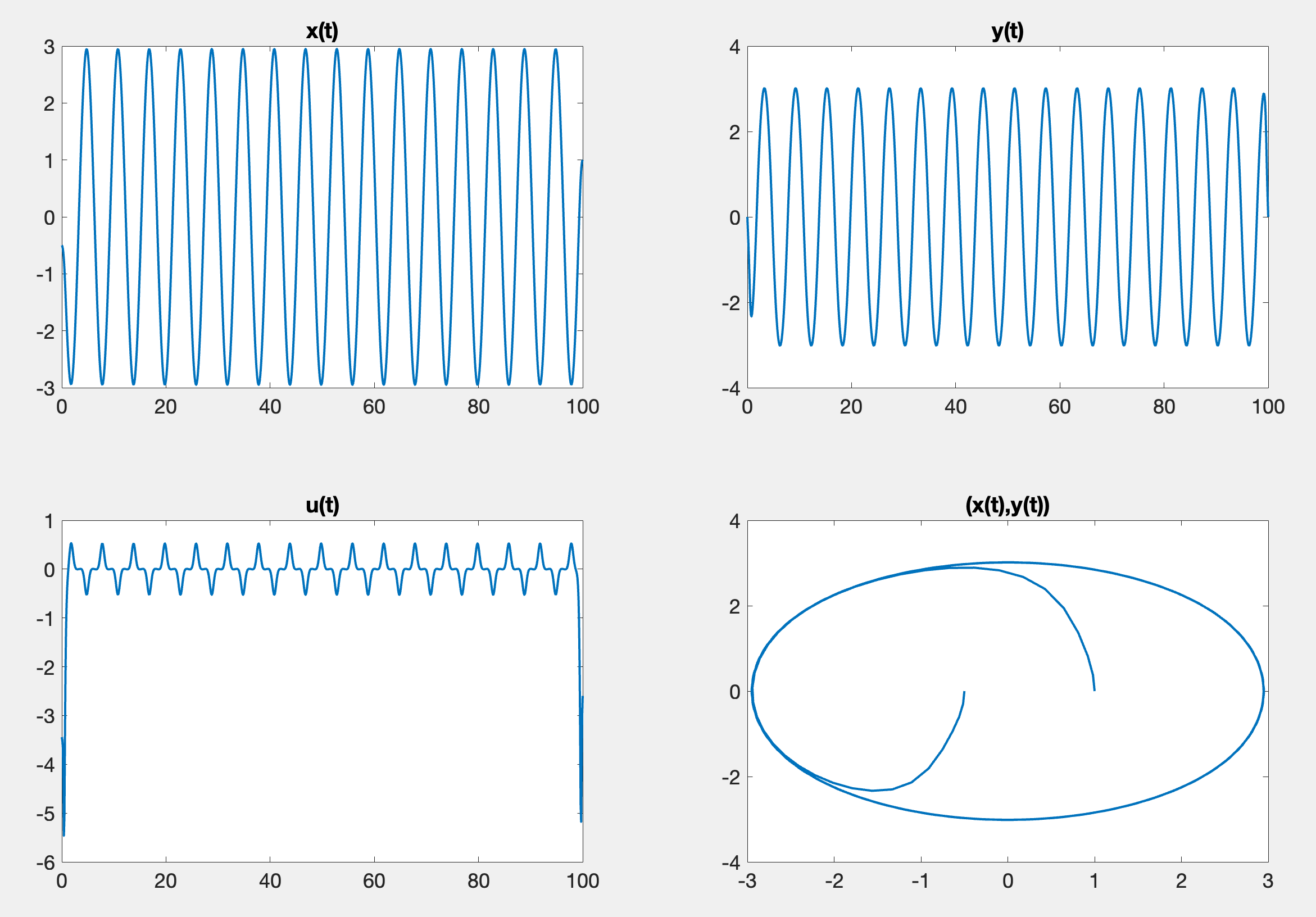}}
\caption{Globally optimal trajectory steering $(x_0,y_0)$ to $(x_1,y_1)$ in time $T=100$, of cost $\simeq 45.6722$.}\label{fig_simu_circle2}
\end{figure}

Numerical simulations illustrating this result are given on Figures \ref{fig_simu_circle1} and \ref{fig_simu_circle2}, with
$$
(x_0,y_0)=(-0.5,0), \qquad (x_1,y_1)=(1,0), \qquad T=100, \qquad R=3.
$$

The computations were performed using AMPL~\cite{AMPL} combined with IpOpt~\cite{IPOPT}, with appropriate initializations (see the comments in Section~\ref{sec_shooting}).

Figure \ref{fig_simu_circle1} shows the locally optimal trajectory steering $(x_0,y_0)$ to $(x_1,y_1)$ in time $T$, which enjoys the local exponential turnpike property near the turnpike $(\bar x,\bar y,\bar u)=(1/2,0,1/2)$ (in red).

Figure \ref{fig_simu_circle2} shows the globally optimal trajectory steering $(x_0,y_0)$ to $(x_1,y_1)$ in time $T$, oscillating closely around the circle $\mathcal{C}$.

\paragraph{More chaotic examples.}
The previous example actually enters in the framework of \emph{periodic turnpike} studied in \cite{TrelatZengZhang, TrelatZhang_MCSS2018, TrelatZhangZuazua_SICON2018}: the circle $\mathcal{C}$ is the turnpike set and the turnpike consists of a periodic trajectory. But there exist much more chaotic examples, as reported in \cite{TrelatZuazua_ongoing}. For example, consider the control system $\dot x(t) = X(x(t)) + b u(t)
$ in $\R^n$, where $b\in\R^n\setminus\{0\}$ and $X$ is a smooth vector field in $\R^n$. We assume that, near the origin, $X(x)=Ax$ with $A$ a $n\times n$ matrix such that $(A,b)$ satisfies the Kalman condition. Far from the origin, we assume that, in some corona $\{x\in\R^n\ \mid\ R_1<\Vert x\Vert<R_2\}$, there are some ergodic trajectories, i.e., filling almost all the corona in infinite time, but no periodic trajectories. %Such chaotic motions do exist and are well known in KAM theory. We take $X$ anyway so that any point of $\R^n$ can be steered to the corona in finite time, in forward or backward time.
Now, take $f^0(x,u)=g(x)+u^2$ with a smooth nonnegative function $g$ that is quadratic near the origin and that is identically zero inside the corona. 
Choose $X$ and $g$ so that the static optimization problem has either a unique local (and thus global) minimizer, or has several local minimizers, with a positive value of $f^0$ at any local minimizer. Then, by construction, we have no global turnpike at any of those minimizers. Any globally optimal trajectory in large time spends most of its time in the corona and has a chaotic behavior.
%
%Hence, the answers to the three questions raised at the beginning of Section \ref{sec_global} are negative.

\paragraph{A meaningful example in dimension one.}
We consider the 1D optimal control problem
\begin{equation}\label{pb_Enrique_1D}
\begin{split}
& \dot x(t) = 4x(t) + \alpha x(t)^3 + u(t) , \\
& x(0)=x_0, \quad x(T)=x_1, \\
& \min \int_0^T \big[ (x(t)-1)^2 + u(t)^2 \big]\, dt ,
\end{split}
\end{equation}
for some $\alpha\in\R$, with fixed terminal conditions $x_0\in\R$ and $x_1\in\R$.
There is no constraint on the (scalar) control.
It follows from known general results (see \cite{Cesari} or \cite[Theorem 2.9]{Trelat_SB}) that, for every $T>0$, there exists at least one optimal solution $(x_{_T}(\cdot),u_{_T}(\cdot))$ of \eqref{pb_Enrique_1D} on $[0,T]$.
By the Pontryagin maximum principle, there exists $\lambda_{_T}(\cdot):[0,T]\rightarrow\R$ such that $u_{_T}(\cdot)=\lambda_{_T}(\cdot)$ and
\begin{equation*}%\label{pb_Enrique_extremal_system_1D}
\begin{split}
\dot x_{_T} &= 4x_{_T}+\alpha x_{_T}^3+\lambda_{_T} , \\
\dot\lambda_{_T} &= -4\lambda_{_T} - 3\alpha x_{_T}^2\lambda_{_T} + x_{_T}-1 .
\end{split}
\end{equation*}
We choose $\alpha \simeq -1.0225539756$, so that the static optimization problem has two global minimizers 
$\bar x_1\simeq 0.059$ and $\bar x_2\simeq 1.961$, and one local minimizer $\bar x_{loc} \simeq -1.925$.
%
%\begin{figure}[H]
%\centerline{\includegraphics[width=7cm]{fig_poly.png}}
%\caption{Graph of $x\mapsto(x-x_d)^2+(ax+\alpha x^3)^2$ for $x_d = 1$, $a = 4$, $\alpha \simeq -1.0225539756$.}\label{fig_poly}
%\end{figure}
%
%\noindent The corresponding values for the control and Lagrange multiplier are, respectively,
%$$
%\bar u_{loc}=\bar\lambda_{loc}\simeq 0.396315259,
%$$
%and
%$$
%\bar u_1=\bar\lambda_1\simeq -0.2358751835,
%\qquad
%\bar u_2=\bar\lambda_2\simeq -0.1264923589.
%$$
%Assumptions \ref{H1_U_W} and \ref{H1_Kalman} of Theorem \ref{thm_TZ} are checked numerically for the three minimizers. Assumption \ref{R_nonsing} is obviously satisfied since $R(x,x'))=(x-x_0,x'-x_1)$.

\begin{figure}[H]
\centerline{\includegraphics[width=8cm]{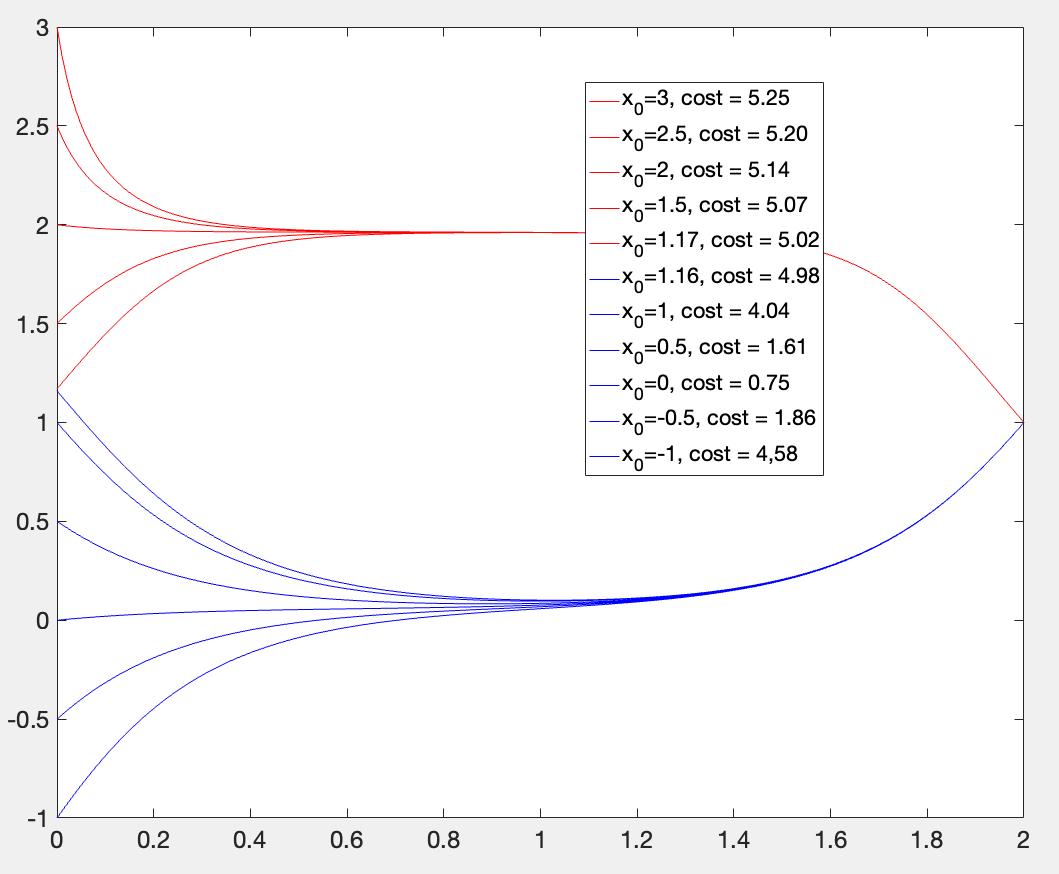}}
\caption{Globally optimal trajectory steering $x_0$ to $x_1=1$ in time $T=2$, for several values of $x_0$.}\label{fig_simu1D_x0}
\end{figure}

\medskip
On Figure \ref{fig_simu1D_x0}, we report numerical simulations done with the fixed final state $x_1=1$ and with $T=2$. We make vary $x_0$ and compute the corresponding globally optimal trajectory.
%The value of $T$ is not large but we already see very well the turnpike phenomenon. Actually, for $T$ larger we obtain the same kind of results but with trajectories approaching the turnpike much quicker.

For $x_0$ less than the critical value that is between $1.16$ and $1.17$, the optimal trajectory enjoys a (global) turnpike property at $\bar x_1$. For $x_0$ beyong this critical value, the turnpike is at $\bar x_2$.

\begin{figure}[H]
%\begin{center}
%\resizebox{10cm}{!}{\input fig_simu1D_partition.pdf_t}
%\end{center}
\centerline{\includegraphics[width=7.5cm]{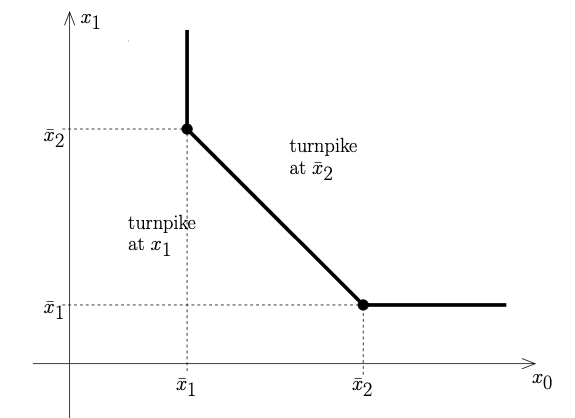}}
\caption{Global turnpike property, depending on $(x_0,x_1)$.}\label{fig_simu1D_partition}
\end{figure}

On Figure \ref{fig_simu1D_partition}, we report the result of a series of numerical simulations, showing how the global turnpike property depends on $x_0$ and $x_1$, in the plane $(x_0,x_1)$. Numerically, this has been obtained by computing the globally optimal trajectory for all $(x_0,x_1)$ ranging in a sufficiently fine grid. 
Interestingly, we note that the bifurcation between a turnpike property at $\bar x_1$ and a turnpike property at $\bar x_2$ seems to be along a polygonal line, made of two (vertical and horizontal) half-lines and of a line segment connecting the two points $(\bar x_1,\bar x_2)$ and $(\bar x_2,\bar x_1)$ in the plane $(x_0,x_1)$.
We do not know how to explain this numerical observation.

\begin{figure}[H]
\centerline{\includegraphics[width=12cm]{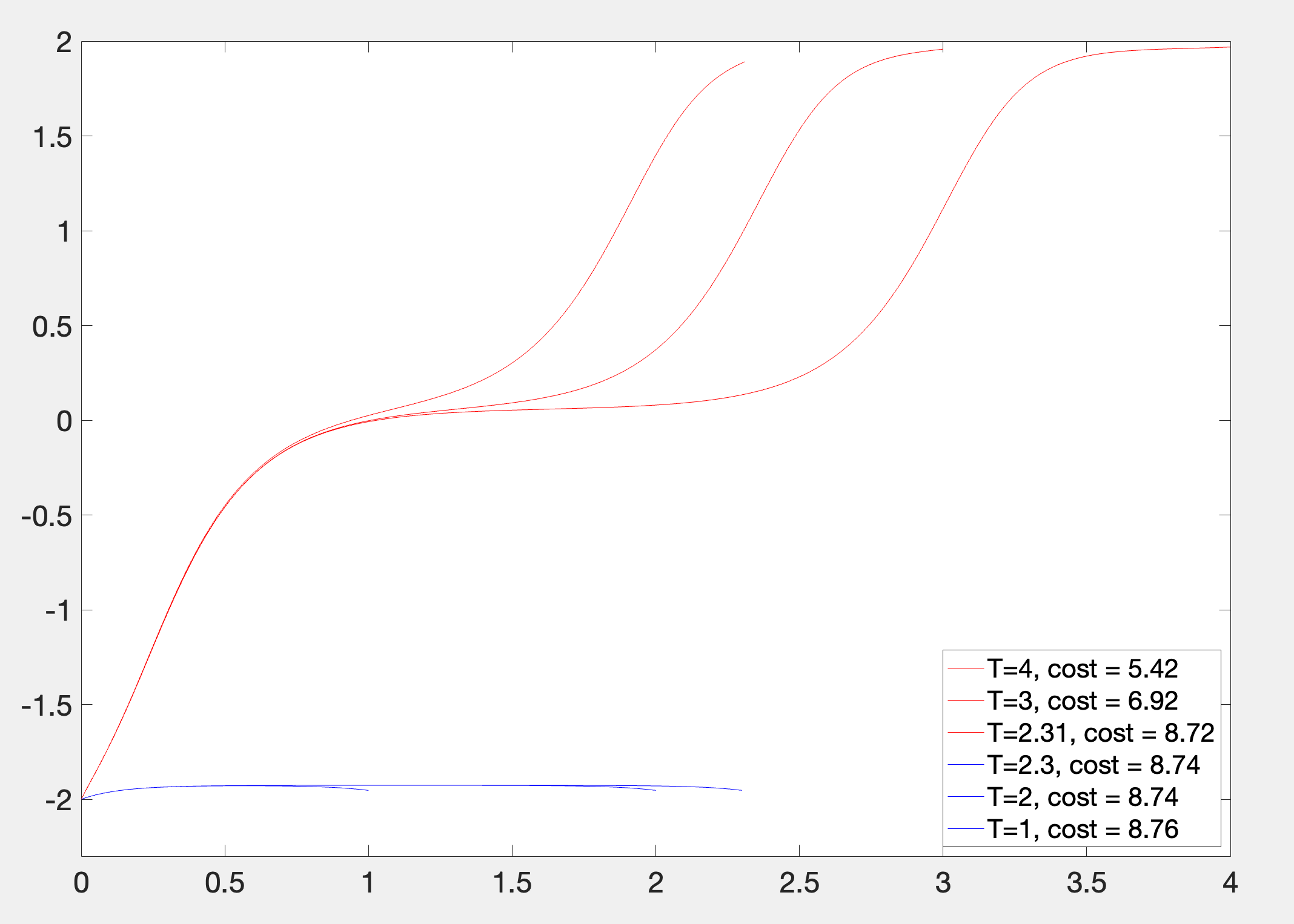}}
\caption{Globally optimal trajectory starting at $x_0=-2$ (with $x(T)$ free), for various values of $T$.}\label{fig_simu1D_T}
\end{figure}

To complete, we illustrate the ``local to global" turnpike phenomenon, depending on the values of $T$. 
On Figure \ref{fig_simu1D_T}, we compute the optimal trajectory of \eqref{pb_Enrique_1D}, with $x_0=2$ and where now $x(T)$ is let free, for various values of $T$. For $T$ small enough (less than $2.3$ on the figure), the optimal trajectory enjoys a turnpike property at the local turnpike $\bar x_{loc}$. For $T$ large enough (more than $2.31$ on the figure), the optimal trajectory enjoys a turnpike property at the global turnpike $\bar x_1$; this is expected because for $T$ large enough a global turnpike becomes less costly.

\paragraph{Open problem.}
It may be interesting to study the following generalization in $\R^n$:
\begin{equation*}
\begin{split}
& \dot x(t) = Ax(t) + \alpha\Vert x(t)\Vert^2x(t) + b\, u(t) ,\\
& x(0)=x_0, \quad x(T)=x_1, \\
& \min \int_0^T \big( \Vert x(t)-x_d\Vert^2 + u(t)^2 \big)\, dt ,
\end{split}
\end{equation*}
where $T>0$, $A$ is a $n\times n$ matrix, $\alpha\in\R$, $b\in\R^n$, with $(A,b)$ satisfying the Kalman condition, $x(t)\in\R^n$ is the state and $u(t)\in\R$ is the control. The initial state $x_0\in\R^n$ and the final state $x_1\in\R^n$ are fixed.
Here, $\Vert\ \Vert$ is the Euclidean norm. 

An interesting situation may be when the matrix $A$ is the discretized Laplacian. In this case, the above system is an approximation of the heat equation.

\section{Generalizations and applications: turnpike everywhere}\label{sec_generalizations}
We have established the exponential turnpike property for finite-dimensional LQ optimal control problems in an elementary way in Section \ref{sec_LQ}, and we have shown a (local) generalization to finite-dimensional nonlinear optimal control problems in Section \ref{sec_nonlinear}.
The turnpike phenomenon withstands a number of other generalizations that we describe hereafter, in infinite dimensions, in the stochastic setting, in mean field games, in optimal design, in situations where the turnpike is an evolving curve and not necessarily a steady-state, etc.

\subsection{Turnpike in infinite dimension}
In infinite dimension, compared with the LQ setting considered in Section \ref{sec_LQ}, $A:D(A)\rightarrow X$ is now a linear operator on the Banach space $X$, generating a strongly continuous semigroup, and $B$ is a control operator. 
This framework involves linear controlled partial differential equations (PDEs), like heat or wave equations with internal or boundary controls. Before mentioning the existing literature, it is worth noting that, for a linear controlled PDE, the control operator $B$ is bounded when the control is applied internally, and it is unbounded when the control is applied at the boundary of the domain.

In this infinite-dimensional LQ context, the turnpike property has been first established in \cite{PorrettaZuazua, PZ2} for a bounded control operator $B$, thanks to Riccati theory. 
 It has been generalized to the case of unbounded control operators in several papers. In \cite{Gugat_2019, GugatHante_2019, GugatTrelatZuazua_SCL2016}, the exponential turnpike property was established for some classes of hyperbolic boundary control problems involving, for instance, wave equations with boundary control.
In \cite{TrelatZhangZuazua_SICON2018}, the exponential turnpike property was established for general abstract LQ optimal control problems, using Riccati theory, either for bounded control operators or for analytic semigroups and unbounded control operators assumed to be admissible (see \cite{Trelat_SB, TucsnakWeiss}), and under exponential stabilizability and detectability assumptions. 
This result was generalized in \cite{GruneGuglielmi_MCRF2021, GruneSchallerSchiela_JDE2020}, and in \cite{GruneSchallerSchiela_COCV2021, Pighin} for some semilinear PDEs.
We also refer to \cite{BreitenPfeiffer, Gugat_MCSS2021, GugatHeitschHenrion, GugatHerty_COCV2023, GugatLazar_2023, WarmaZamorano, Zamorano} for related results.

\subsection{Turnpike in mean field games}%\label{sec_MF}
The authors of \cite{Cardaliaguet2012, Cardaliaguet2013, LasryLions} investigated the large-time behavior of solutions of mean field games (MFG)
\begin{equation*}
\left\{
\begin{split}
&-\partial_t v - \triangle v + H(x,Dv) = F(x,m) \qquad \textrm{in }(0,T)\times\mathbb{T}^d , \\ 
& \partial_t m - \triangle m - \mathrm{div}(m\, H_p(x,Du)) = 0 ,  \\  % 
& m(0)=m_0,\quad v(T)=G(\cdot,m(T)) .
\end{split}\right.
\end{equation*}
In this system, the first equation is the Bellman equation for the agents' value function $v(t,\cdot)$.
The second equation is the Kolmogorov-Fokker-Planck equation for the distribution of agents, $m(t,\cdot)$ is the probability density of the state of players at time $t$.
They established large-time inequalities that, in retrospect, correspond precisely to turnpike inequalities, although this interpretation was not made explicit at the time.

Since then the turnpike phenomenon has been well identified for MFG, see \cite{CirantPorretta_COCV2021}.
Mainly, they establish that, if $H$ is uniformly convex in $p$ and if $F$ and $G$ are monotone (even mildly non-monotone) then
$$
\Vert m-\bar m\Vert_{L^\infty} + \Vert Dv-D\bar v\Vert_{L^\infty} \leq C \big( e^{-\nu t} + e^{-\nu (T-t)} \big) \qquad \forall t\in[0,T] ,
$$
which is an exponential turnpike inequality, where the ``turnpike" is the steady-state solution of
\begin{equation*}
\begin{split}
& \bar\lambda - \triangle\bar v + H(x,D\bar v) = F(x,\bar m) \\ 
& -\triangle\bar m - \mathrm{div}(\bar m\, H_p(x,D\bar v)) = 0
\end{split}
\end{equation*}
Moreover
$$
v(t,x) = \bar\lambda(T-t) + \bar v(x) + \theta + \mathrm{o}(1) \qquad \textrm{as}\ T\rightarrow+\infty
$$
uniformly on compacts, where $\theta$ is the unique ergodicity constant of some ergodic stationary problem (this is related to some facts and references mentioned in Section \ref{sec_KAM}).

\subsection{Turnpike, randomness and mean field}
Averaged turnpike properties have been established in~\cite{HernandezLecarosZamorano, HernandezZuazua} in the presence of random coefficients. We also highlight the recent contributions~\cite{SunWangYong, SunYong_JDE2024}, which address the exponential turnpike property in a stochastic setting.

For mean-field dynamics, such as Vlasov equations, we refer to~\cite{Averboukh, GugatHertySegala_2024}, where turnpike inequalities are derived both for sequences of particle systems and their mean-field limits. In the deterministic mean-field optimal control setting, exponential turnpike theorems have also recently been established in~\cite{BonnetWeillColomboShishmintsevTrelat_2026}.

We simply note here that this is a rich and rapidly evolving area of research, which certainly deserves further exploration.

\subsection{More general turnpikes}\label{sec_moregeneral}
Until now we have reported on the turnpike phenomenon where the turnpike consists of a single point. But the same qualitative behaviour can be found in other, more general, settings.

\paragraph{Several turnpike points.}
On the one part, there may exist several turnpike points, and in this case, interestingly, there may exist a competition between the several turnpikes for optimal trajectories having different initial states (see \cite{Rapaport, Rapaport2}, see \cite[Appendix]{Trelat_MCSS2023}, see also  Section \ref{sec_global} above).

\paragraph{Periodic turnpike.}
On the other part, for some classes of optimal control problems with periodic systems, the turnpike phenomenon may occur around a periodic trajectory, which is itself characterized as being the optimal solution of an appropriate periodic optimal control problem (see \cite{Samuelson1976, TrelatZengZhang, TrelatZhangZuazua_SICON2018, ZanonGruneDiehl, Zaslavski_2006, Zaslavski_2015}).
For instance, the exponential periodic turnpike inequality of \cite{TrelatZhangZuazua_SICON2018} was established by using a periodic version of the Riccati theory. 
There are many examples of periodic turnpike in biology, because in this field periodic phenomena are ubiquitous. 

\paragraph{Partial turnpike.}
There are situations in which the optimal control problem does not admit a full equilibrium but rather a \emph{partial} equilibrium, as studied in~\cite{FaulwasserFlasskampOberBlobaumSchallerWorthmann_MCSS2022, Trelat_MCSS2023, FlasskampMaslovskayaOberBlobaumWembe_MCSS2025}. In such cases, the turnpike set may take the form of a curve. A basic example, which motivated the study in~\cite{Trelat_MCSS2023}, is a runner model where the objective is to complete a $1500\,\mathrm{m}$ race while optimizing a cost that typically balances time and energy expenditure.

\begin{figure}[h]
\centerline{\includegraphics[width=6cm]{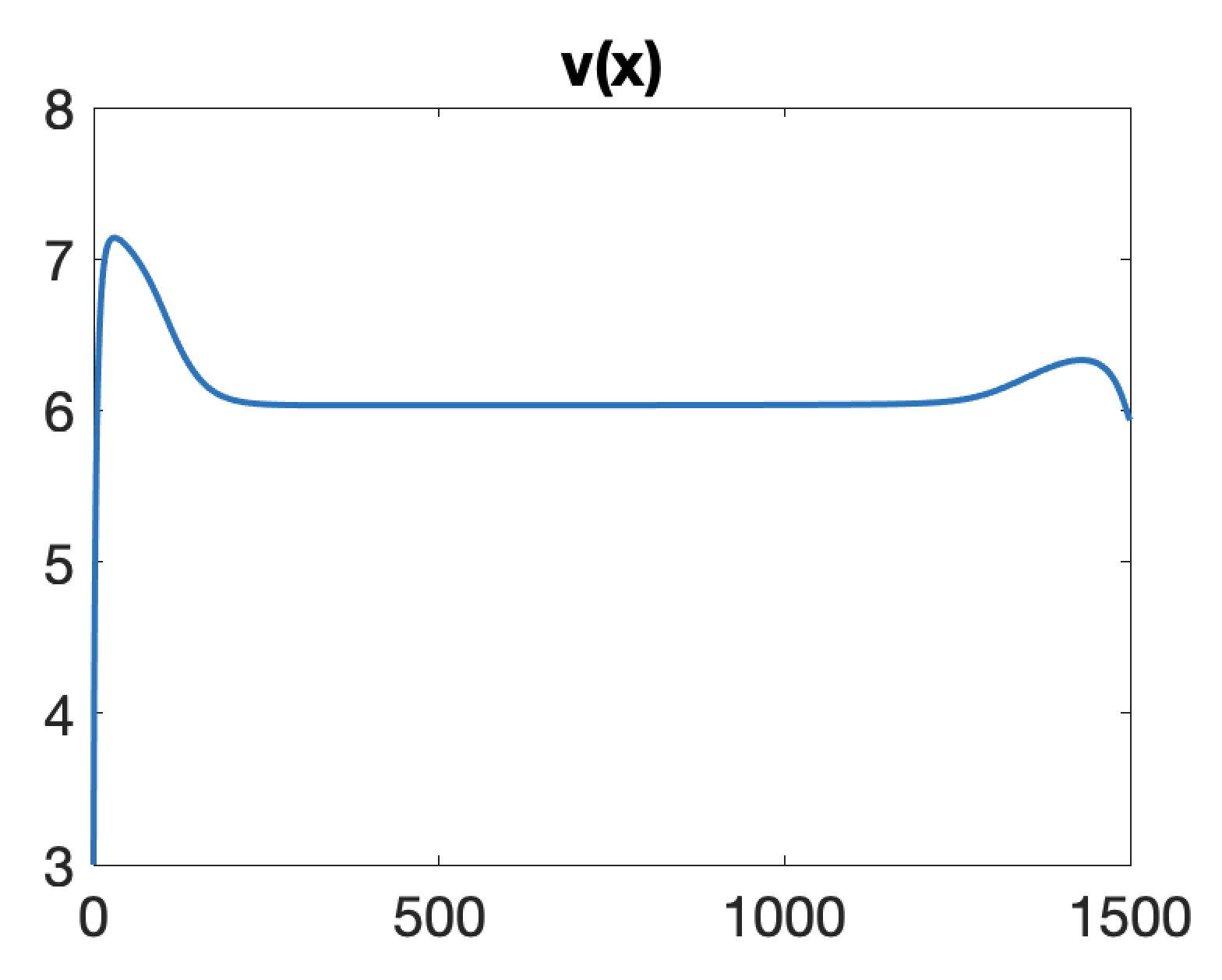}}
\caption{Speed of a runner.}\label{fig_runner}
\end{figure}

In this context, except at the start and the finish of the race, the runner's speed remains essentially constant (see Figure~\ref{fig_runner}). However, the turnpike is not a full equilibrium since the runner's position evolves continuously over time. This application has been further explored in~\cite{AftalionTrelat_RSOS, AftalionTrelat_JOMB}, where the turnpike property plays a crucial role in ensuring the convergence of the numerical algorithms employed (see also Section~\ref{sec_shooting} for further details).

\paragraph{General turnpike set.}More generally, the \emph{turnpike set} can take a very general form. The analysis carried out in~\cite{TrelatZhang_MCSS2018} provides a characterization of the value function near the turnpike set within a broad framework, leveraging dissipativity arguments.

The turnpike set can also be interpreted from the perspective of dynamical systems theory, as explored in~\cite{SakamotoZuazua}. Near a singular point of a dynamical system, it is well known that invariant manifolds arise -- namely, the stable, unstable, and center manifolds. The stable and unstable manifolds precisely capture the hyperbolicity phenomenon, which, as we have shown, can be exploited to establish the exponential turnpike property.

In this framework, the center manifold naturally corresponds to the turnpike set. When the center manifold reduces to a single point, the turnpike is likewise a point. However, when the center manifold is nontrivial, it possesses its own intrinsic dynamics -- potentially periodic (yielding a periodic turnpike) or even chaotic -- referred to as the \emph{turnpike dynamics}. The flows along the stable and unstable manifolds then govern the hyperbolic behavior surrounding the turnpike dynamics of the center manifold.

From this perspective, many behaviors or phenomena across various fields and applications can be reinterpreted as manifestations of the turnpike property. A typical example arises in space applications, where the optimal trajectory of a spacecraft often consists of an initial brief thrust phase, followed by a long coasting period with no thrust (such as a ballistic phase), and a final brief thrust to reach the target. This structure is characteristic of low-thrust orbit transfers or missions involving Lagrange points (see, e.g.,~\cite{ChupinHaberkornTrelat_M2AN2017, ChupinHaberkornTrelat_JGCD2018}), where the turnpike phenomenon has been explicitly identified and exploited.

We also refer to~\cite{KLMR} for examples of nontrivial center manifolds in celestial mechanics, including cases exhibiting chaotic behavior. In such contexts, the presence of turnpike dynamics could likely be identified, particularly if future space missions were to exploit these chaotic central manifolds.

\paragraph{Polynomial turnpike.} When first-order optimality conditions, such as the Pontryagin Maximum Principle, can be applied, most results in the literature establish an exponential turnpike estimate. In the recent work~\cite{HanZuazua_2022}, an integral turnpike estimate with polynomial rate has been derived for a linear hyperbolic system with quadratic cost, under weak controllability and observability assumptions. More recently,~\cite{TrelatZuyev} proves a pointwise polynomial turnpike property around a steady state for a particular class of infinite-dimensional optimal control problems, including models such as rotating flexible beams.

\subsection{Turnpike and dissipativity}\label{sec_dissipativity}
The concept of dissipativity was introduced by Willems in the 1970s (see~\cite{Willems}). The family of optimal control problems \eqref{syst}--\eqref{syst_Omega}--\eqref{mincost}, indexed by $T > 0$, is said to be \emph{strictly dissipative}\footnote{We refer to \emph{nonstrict} dissipativity when $\alpha = 0$ in \eqref{dissip}.} at the optimal static point $(\bar x, \bar u)$ with respect to the \emph{supply rate function}
\begin{equation*}%\label{supply}
w(x, u) = f^0(x, u) - f^0(\bar x, \bar u)
\end{equation*}
if there exists a locally bounded function $S: \R^n \rightarrow \R$, called the \emph{storage function}, such that the following \emph{strict dissipativity inequality} holds:
\begin{equation}\label{dissip}
S(x(t_0)) + \int_{t_0}^{t_1} w(x(t), u(t))\, dt \geq S(x(t_1)) + \int_{t_0}^{t_1} \alpha\left(\Vert x(t) - \bar x \Vert, \Vert u(t) - \bar u \Vert\right) dt
\end{equation}
for all $t_0 < t_1$ and for any admissible pair $(x(\cdot), u(\cdot))$ solving \eqref{syst} and \eqref{syst_Omega}, where $\alpha$ is a function of class $\mathcal{K}$, i.e., $\alpha: [0, +\infty) \rightarrow [0, +\infty)$ is continuous, increasing, and satisfies $\alpha(0) = 0$. Further discussions on dissipativity can be found in~\cite{Brogliato, TrelatZhang_MCSS2018, Willems}.

While the above definition is global, local versions of dissipativity can naturally be considered on subsets of $\R^n$.

Dissipativity has been widely used in~\cite{DammGruneStielerWorthmann_SICON2014, FaulwasserKellett_A2021, FaulwasserKordaJonesBonvin_A2017, GruneGuglielmi_SICON2018, GruneGuglielmi_MCRF2021, GruneKellettWeller, GruneMuller_SCL2016, GruneGuglielmi_SICON2018} to establish turnpike properties.

In essence, for linear-quadratic (LQ) problems, dissipativity is equivalent to the turnpike property. Generalizations to more degenerate settings have been obtained in~\cite{GuglielmiLi_MCSS2024, GugatLazar_A2023, LuGuglielmi}.

More generally, it is shown in~\cite{TrelatZhang_MCSS2018} (both in finite and infinite dimensions) that
\[
\text{strict strong duality} \ \Rightarrow \ \text{strict dissipativity} \ \Rightarrow \ \text{measure-turnpike property}.
\]
We now provide a more precise discussion of these implications.

The static problem \eqref{staticpb} is said to have the \emph{strong duality property} if $(\bar x,\bar u)$  minimizes the \emph{Lagrangian function} $L(\cdot,\cdot,\bar \lambda)$ defined by $L(x,u,\bar \lambda)= f^0(x,u)-\langle \bar \lambda, f(x,u)\rangle$ (the sign minus is due to the fact that we took $\lambda^0=-1$ for the Lagrange multiplier associated with the cost). This notion is well known in classical optimization theory, in relation with primal and dual problems. It is proved in \cite{TrelatZhang_MCSS2018} that the (strict) strong duality property implies the (strict) dissipativity property, with the linear storage function $S(x)=\langle \bar \lambda,x\rangle$.

Strong duality is an infinitesimal version of the dissipativity inequality, at least when the storage function is continuously differentiable (see also \cite{FaulwasserKordaJonesBonvin_A2017, FaulwasserGrune_2022}). See also Section \ref{sec_KAM} for another interesting relationship exploiting this infinitesimal version.

It is also proved in \cite{TrelatZhang_MCSS2018} that the strict dissipativity property implies the measure-turnpike property, i.e., for every $\varepsilon>0$, there exists $\Lambda(\varepsilon)>0$ such that
$$
| \{ t\in[0,T]\, \mid\, \|(x_{_T}(t)-\bar x,u_{_T}(t)-\bar u)\| > \varepsilon \} |\leq \Lambda(\varepsilon),\,\forall T>0,
$$
where the above notation $\vert\cdot \vert$ stands for the Lebesgue measure. We refer the reader to \cite{CarlsonHaurieLeizarowitz_book1991, FaulwasserKordaJonesBonvin_A2017, Zaslavski_2006} (and references therein) for similar definitions.

The measure-turnpike property quantifies the measure of the set of times during which the optimal trajectory and control remain outside an $\varepsilon$-neighborhood of the turnpike. Naturally, this property is much weaker than the exponential turnpike property. Moreover, it is worth noting that the measure-turnpike estimate involves only the state and control variables, while the costate does not appear explicitly in the inequality.

\medskip

As mentioned at the end of Section~\ref{sec_nonlinear}, dissipativity has often been employed in the literature as an alternative to linearization or Riccati theory for establishing the turnpike property, particularly in situations where these classical techniques fail to apply -- such as when the optimal control problem involves state and/or control constraints that are saturated along the turnpike.

We end this section by mentioning another use of (global) dissipativity, as a tool to infer global results from local ones. Theorem \ref{thm_TZ} in Section \ref{sec_nonlinear}, which is the main result of \cite{TrelatZuazua_JDE2015} states an exponential turnpike property that is of a \emph{local} nature, as already discussed. %, because of the smallness condition \eqref{smallness}. 

As noted at the end of Section~\ref{sec_nonlinear} and demonstrated in~\cite{Trelat_MCSS2023}, strengthening Theorem~\ref{thm_TZ} by assuming global dissipativity -- along with uniqueness of the minimizers for both the optimal control problem and the static problem -- leads to a \emph{global} exponential turnpike property. Although the uniqueness assumptions are generic, the global dissipativity assumption is significantly stronger and, in general, challenging to verify.

These difficulties, along with the global aspects of the turnpike property, have been illustrated in Section~\ref{sec_global}.

We also mention the article \cite{EsteveGeshkovskiPighinZuazua_2022}, where bootstrap arguments have been used to pass from local to global, in some specific situations.

\subsection{Turnpike and weak KAM}\label{sec_KAM}
Related to strong duality, the infinitesimal form of the non-strict dissipativity inequality \eqref{dissip}, with $\alpha=0$, is the Hamilton-Jacobi inequality 
$$
H_1(x,\nabla S(x))\leq -f^0(\bar x,\bar u) 
$$
where $H_1$ is the maximized normal Hamiltonian (indeed, divide by $t_1-t_0>0$ and take the limit $t_1-t_0\rightarrow 0$). The existence of $C^1$ solutions is therefore related to \emph{weak KAM theory} (see \cite{Fathi}). In this context, the singleton $\{(\bar x,\bar u)\}$ is the \emph{Aubry set} and $f^0(\bar x,\bar u)$ is the \emph{Ma{\~{n}}\'e critical value}. 

This fact has been noted in \cite{Trelat_MCSS2023} but not yet fully exploited, up to our knowledge.

Accordingly, the turnpike theory has a close relationship with the concept of \emph{ergodicity} for Hamilton-Jacobi-Bellman equations, which is a property used to describe the large-time behavior of the solutions of HJB equations in relation to the solution of a stationary equation (see \cite{Arisawa, BarlesSouganidis, BuckdahnQuincampoixRenault, LiQuincampoixRenault, QuincampoixRenault}).
We also refer to \cite{ATZ_MCRF2024, ATZ_SCL2024, EKPZ} where a two-terms expansion of the value function $V$ has been derived, yielding an asymptotic expansion of the form
\begin{equation}\label{asympt_exp}
V(T,x_0,x_1) = T \, f^0(\bar x,\bar u) + V_f(x_0) + V_b(x_1) + \mathrm{o}(1)
\end{equation}
as $T\rightarrow+\infty$, where
\begin{equation*}
\begin{split}
& V_f(x_0) = \min\int_0^{+\infty} \left( f^0(x,u)-f^0(\bar x,\bar u) \right) dt , \\
& \dot x=f(x,u),\quad x(0)=x_0,
\end{split}
\end{equation*}
is the optimal cost for stabilizing the forward dynamics $f$ from $x_0$ to $\bar x$, and
\begin{equation*}
\begin{split}
& V_b(x_1) = \min\int_0^{+\infty} \left( f^0(x,u)-f^0(\bar x,\bar u) \right) dt, \\
& \dot x=-f(x,u),\quad x(0)=x_1,
\end{split}
\end{equation*}
is the optimal cost for stabilizing the backward dynamics $-f$ from $x_1$ to $\bar x$.
In \eqref{asympt_exp}, the dominating term $T \, f^0(\bar x,\bar u)$ corresponds of course to the cost of remaining at the turnpike (which is the Aubry set, as said above), and the next term in the expansion is the sum of optimal stabilization costs. The expansion \eqref{asympt_exp} can be written more explicity in the LQ case (see \cite{ATZ_MCRF2024}, where the result holds also in infinite dimension).

\subsection{Turnpike in shape optimal design}
The concept of shape turnpike has been introduced in \cite{LanceTrelatZuazua_SCL2020, LanceTrelatZuazua_2022} (see also \cite{AM} for a related issue). Let us explain the concept on a simple model.

Let $\Omega\subset\R^n$ be a bounded domain, let $y_0, y_d\in L^2(\Omega)$, and let $L\in(0,1)$. Given any $T>0$, we consider the following dynamical optimization problem for the Dirichlet heat equation on $\Omega$ with source term:
\begin{equation*}
\begin{split}
&\partial_t y = \triangle y + \chi_{\omega(t)}, \\
&y_{\vert\partial\Omega}=0, \\
&y(0)=y_0, \\
&\min_{\omega(\cdot)} \int_0^T \Vert y(t)-y_d\Vert_{L^2(\Omega)}^2\, dt  ,
\end{split}
\end{equation*}
where the minimum is taken over all measurable subsets $\omega(t)\subset\Omega$ of measure $\vert\omega(t)\vert =L\vert\Omega\vert$.
Here, we have set 
$$
\chi_\omega(x) = \left\{\begin{array}{ll} 1 & \textrm{if}\ x\in\omega, \\ 0 & \textrm{otherwise.} \end{array}\right.
$$
In the above problem, the control is the subdomain $\omega(t)$, that one can choose at any instant of time in order to minimize the integral criterion. This situation falls outside the scope of standard PDE control, which is why we refer to it as \emph{shape control}.

It is natural to expect a turnpike phenomenon for large $T$. Accordingly, we consider the static shape optimization problem
\begin{equation*}
\begin{split}
& \triangle y + \chi_\omega = 0, \\
& y_{\vert\partial\Omega}=0 \\
& \min_{\omega} \Vert y-y_d\Vert_{L^2(\Omega)}^2 ,\\
\end{split}
\end{equation*}
over all possible measurable subsets $\omega\subset\Omega$ of measure $\vert\omega\vert =L\vert\Omega\vert$.

It is proved in \cite{LanceTrelatZuazua_SCL2020} that, if either $y_d$ is convex or $y_d\leq 0$ or $y_d\geq \bar y^{1}$ (with $\triangle\bar y^{1}+1=0$, $\bar y^{1}_{\vert\partial\Omega}=0$) then the dynamical shape optimization problem has a unique solution $t\mapsto \omega(t)$, the static shape optimization problem has a unique solution $\bar\omega$, and there exists $M>0$ such that
$$
\int_0^T \left( \Vert y(t)-\bar{y} \Vert_{L^2(\Omega)}^2 + \Vert \lambda(t)-\bar\lambda \Vert_{L^2(\Omega)}^2 + \Vert \chi_{\omega(t)}-\chi_{\bar\omega} \Vert_{L^1(\Omega)}\right) dt \leq M \qquad \forall T>0
$$
(\emph{shape turnpike} inequality, which is, here, an \emph{integral turnpike} property), where $\lambda(t)$ and $\bar\lambda$ are the costates / Lagrange multipliers of the respective optimization problems (they are obtained by applying the bathtub principle). As a consequence, we have
$$
\frac{1}{T}\int_0^T y(t)\, dt \underset{T\rightarrow+\infty}{\overset{L^2}{\longrightarrow}}\bar y,
\quad 
\frac{1}{T}\int_0^T \lambda(t)\, dt \underset{T\rightarrow+\infty}{\overset{L^2}{\longrightarrow}}\bar\lambda,
\quad
\frac{1}{T}\int_0^T \chi_{\omega(t)}\, dt \underset{T\rightarrow+\infty}{\overset{L^1}{\longrightarrow}}\chi_{\bar\omega} .
$$ 
It is not known whether an exponential turnpike property can be obtained. The proof uses dissipativity considerations.

\subsection{Turnpike in aerospace and aeronautics}
As alluded to in Remark \ref{rem4}, a well known and quite evident example of an ``exact" turnpike is provided by the orbit transfer problem, where the minimal-consumption trajectories typically consist of three successive arcs: a first, short-time, thrust arc (along which the thrust is maximal), then a long ``coast arc" (also called ``balistic arc", along which there is no thrust), and finally a short-time thrust arc again.
Such strategies are widely used by space agencies to perform orbit transfers after having launched an engine. Sometimes, strategies may involve a larger number of thrust arcs. Anyway, in this picture, the middle coast (balistic) arc is viewed as an exact turnpike; it is exact because this arc corresponds to a free-motion Kepler orbit, solution of the classical Newton dynamics.

A variant of this problem, with the Kepler model but where the authors have considered a quadratic cost functional, has been studied in \cite{FlasskampMaslovskayaOberBlobaumWembe_MCSS2025}; accordingly, the turnpikes are free-motion Kepler orbits. 

\medskip

It is interesting to note that turnpike phenomena have been observed in shape optimization problems in aeronautics. Indeed, in many aeronautics problems the optimal design is performed by engineers based on stationary models (see \cite{Jameson, JamesonMartinelli}). This remarkable engineering intuition can actually be rigorously justified by the turnpike phenomenon, which allows one to determine optimal shape designs or optimal materials on the basis of a steady-state model.

For example, in a basic sense, this is similar to how an airplane operates: the pilot is extremely active and dynamic during takeoff and landing, but the cruise phase is essentially stationary.
Again, this clearly reflects the turnpike phenomenon.
Similar features appear in an ubiquitous way in many problems in aeronautics. 

In more general, the turnpike principle can be used to tackle in a more tractable way some complex parameter dependent optimal control problems, or can be used as a model reduction.
In that view, turnpike offers a natural approximation to the long-time problem under consideration, by reducing it to a steady state.
As said above, this ``turnpike model reduction" has been used quite much in aeronautics, particularly in shape optimization problems and parameter identification issues. 

This is the case in \cite{ErsoyFeireislZuazua_JDE2013}, where, in view of optimal parameter identification within the framework of hyperbolic conservation laws, the authors consider 1D forced steady-state scalar conservation laws and perform a sensitivity analysis of the shock location with respect to variations of the forcing term. 

Another interesting and challenging problem, that we would like to mention, in which this paradigm could be used, is the sonic-boom minimization problem (see \cite{AlonsoColonno_2012}).
In this problem, the objective is to develop supersonic aircrafts that are sufficiently quiet to be allowed to fly supersonically over land. The pressure signature created by the aircraft is expected to be such that,
when reaching ground, it can barely be perceived by humans, and it results in admissible disturbances to man-made structures that do not exceed a given threshold of annoyance for the inhabitants. 

In order to be able to design an aircraft with sufficiently low sonic booms, one has to feature how pressure disturbances appear and propagate in the atmosphere, while remaining below an acceptable threshold at the ground. At the same time, from the viewpoint of material optimization, one should determine a fair balance so that the created aircraft be economically and environmentally acceptable. 

This question of control and shape optimization design in long time actually falls under the turnpike theory and could certainly be analyzed at the light of some results that we have given in this paper.

\subsection{Turnpike and numerical optimal control}\label{sec_shooting}
As shown in~\cite{TrelatZuazua_JDE2015}, the exponential turnpike inequality proves particularly useful for initializing numerical methods in optimal control, significantly improving their chances of successful convergence. The underlying idea is that, since the turnpike provides a good approximation of the optimal trajectory, any numerical algorithm should ideally be initialized near the turnpike -- provided it can be computed in advance.

In this sense, the turnpike serves as a tool for complexity reduction, effectively eliminating the explicit dependence on time.
There are, mainly, two different types of numerical methods in optimal control (see \cite{TrelatJOTA} for a review): direct and indirect methods.

Roughly speaking, direct methods consist of discretizing the whole optimal control problem: one chooses a numerical scheme to discretize the dynamics (Euler, Runge-Kutta, etc; or a collocation method) and a numerical scheme to discretize the cost functional (if it is in integral form: rectangle, trapeze, Simpson, etc). All in all, after discretization, the optimal control problem is approximated by a finite-dimensional optimization problem with equality and inequality constraints. This optimization problem can then be solved using classical algorithms, using primal or dual methods (see \cite{Betts}). This so-called direct approach is: first discretize, then optimize (or dualize).

In contrast, indirect methods consist of first optimizing (or dualizing) and then discretizing. This means that we first apply the Pontryagin maximum principle to the optimal control problem, thus reducing it, thanks to the introduction of the costate and taking into account the transversality conditions, to a boundary value problem settled on the extremal system in the variable $z_{_T}$, where $z_{_T}(t)=(x_{_T}(t),\lambda_{_T}(t))$ is the pair (state, costate). This two-point boundary value problem is then in general solved thanks to a Newton method. The resulting algorithm thus combines a Newton method with a method of integration of differential equation: this is the so-called shooting method (see \cite{StoerBulirsch}).

Of course, there are numerous variants of both direct and indirect methods, and some approaches lie between the two and cannot be classified strictly as either. Regardless of the method, initialization is always required, and when the turnpike is known, using it as an initial guess is certainly a sensible strategy.

When the turnpike can be computed, initializing direct methods becomes relatively straightforward. For shooting methods, however, the initialization is less obvious. In the specific case where the turnpike reduces to a single point, a variant of the classical shooting method was introduced and implemented in~\cite{TrelatZuazua_JDE2015}.

In the standard shooting method, one guesses $z_{_T}(0)$ and applies a Newton method to adjust $z_{_T}(T)$ so that the boundary conditions are satisfied. In the variant proposed in~\cite{TrelatZuazua_JDE2015}, the shooting is instead initialized at the midpoint of the trajectory. At $t = T/2$, the optimal $z_T(T/2)$ is exponentially close to $(\bar x, \bar\lambda)$, making it natural to initialize exactly at $(\bar x, \bar\lambda)$. The extremal system is then integrated backward over $[0, T/2]$ to compute $z_T(0)$, and forward over $[T/2, T]$ to compute $z_T(T)$. Finally, the value of $z_T(T/2)$ is adjusted via a Newton method to satisfy the terminal conditions.

This variant of the single shooting method has proven highly efficient. Examples illustrating the effectiveness of turnpike-based initialization can be found in~\cite{AftalionTrelat_JOMB, CaillauDjemaGouzeMaslovskayaPomet, CaillauFerrettiTrelatZidani, Trelat_MCSS2023, TrelatZuazua_JDE2015}.
\begin{example}\label{ex1}
Consider the optimal control problem \eqref{exa} studied in Section \ref{sec_ex1}.
The problem is apparently simple because it seems to be close to a LQ optimal control problem. However the nonlinear term $y^3$ in the second equation is far from being harmless.
According to the above discussion, whether one chooses a direct or an indirect approach to solve numerically the optimal control problem \eqref{exa}, it will always be relevant to initialize the numerical method by using the turnpike solution $(\bar x,\bar u,\bar\lambda)$. Of course, then, in accordance with Theorem \ref{thm_TZ}, the numerical method is expected to converge successfully for initial and final states in the neighborhood of the turnpike under consideration.
The numerical simulations provided on Figure \ref{fig_simu_xcube} were been obtained in two ways. 

First, we implemented a direct approach, by discretizing \eqref{exa} in time with a Crank-Nicolson scheme (implicit trapezoidal). We used the automatic differentiation software AMPL (see \cite{AMPL}) combined with the optimization routine IpOpt (see \cite{IPOPT}). For $T$ sufficiently large, this may fail if we take a quite random initialization, while it always works successfully if we initialize with the turnpike values.

Second, we implemented a single shooting method. We observed that it is impossible to make converge the usual shooting method if $T$ becomes too large, while if we use the variant described above, initialized at the middle with the turnpike values, then convergence is easily obtained, whatever the value of $T$. Here, the impact of the nonlinear term $y^3$ in the dynamics is very significant.
This illustrates the efficiency of this variant of the shooting method.
\end{example}

\subsection{Turnpike in deep learning}
It is well known that deep learning can be formulated as an optimal control problem, where backpropagation corresponds precisely to the backward-in-time costate equation in the Pontryagin Maximum Principle (see~\cite{LeCun}). Consider a neural network with $N+1$ layers and a dataset $\{x^0_k \in \R^n, \ k = 1, \ldots, K\}$ of $K$ input samples. The learning problem can be written as
\begin{equation*}
\begin{split}
& x^{i+1}_k = \sigma(A_i x^i_k), \qquad i = 0, \ldots, N, \quad k = 1, \ldots, K, \\
& \min_{(A_i)_{i=0}^N} \mathrm{loss}(x^N) = \frac{1}{K} \sum_{k=1}^K g_k(x^N_k),
\end{split}
\end{equation*}
where each $A_i$ is the matrix of weights (i.e., the parameters representing interactions in the neural network), and $\sigma$ is a nonlinear activation function (typically a sigmoid or ReLU function).

This formulation clearly defines a discrete-time optimal control problem, where the control variables are the weights $u = (A_i)_{i=0}^N$. 

The residual variant takes the form
\begin{equation*}
\begin{split}
& x^{i+1}_k = x^i + \delta \sigma(A_i x^i_k), \qquad i = 0, \ldots, N, \quad k = 1, \ldots, K, \\
& \min_{(A_i)_{i=0}^N} \mathrm{loss}(x^N) = \frac{1}{K} \sum_{k=1}^K g_k(x^N_k),
\end{split}
\end{equation*}
with $\delta >0$.

When $N$ is large, i. e., in the deep regime, this problem is an approximation of the continuous-time optimal control problem
$$
\dot x_k(t)=\sigma(u(t)x_k(t)),\qquad \min \mathrm{loss}(x(T)),
$$
with $T$ large,
which is a continuous-time optimal control problem in large time, thus expected to enjoy a turnpike property.

This idea has been explored in \cite{EsteveGeshkovskiPighinZuazua_2020} where it has been indeed shown that the turnpike phenomenon occurs. This means that, except at the boundary layers, not much occurs in the network and thus it is not better to use very deep networks to obtain better results.
We refer to the recent work \cite{FaulwasserHempelStreif_2024} for an estimate on appropriate depth bounds.

 \subsection*{Acknowledgements}
This project has received funding from the European Research Council (ERC) under the European Union’s Horizon Europe research and innovation programme (grant agreement No. 101096251-CoDeFeL). This material is based upon work supported by the Air Force Office of Scientific Research under award number FA8655-22-1-7012. EZ was partially supported by the Alexander von Humboldt Professorship program; the European Union’s Horizon Europe MSCA project ModConFlex (HORIZON-MSCA-2021-DN-01, project 101073558); the Transregio 154 Project Mathematical Modeling, Simulation and Optimization Using the Example of Gas  Networks of the DFG; SURE-AI: The Norwegian Centre for Sustainable, Risk-Averse, and Ethical AI, grant 357482, Research Council of Norway; the Madrid Government-UAM Agreement for the Excellence of the University Research Staff in the context of the V PRICIT;  Grant PID2023-146872OB-I00-DyCMaMod of MICIU (Spain); and by the COST Actions "CA24122-Multiscale Stochastics, Patterns, and Analysis of Combinatorial Environments" and "CA24136-Interactions between Control Theory and Machine Learning".

\end{document}